\documentclass[11pt]{amsart}

\usepackage{url}

\usepackage[marginratio=1:1,height=642.4pt,width=391pt,tmargin=86.6pt]{geometry}

\newcommand{\Zbl}[1]{Zbl #1}
\newcommand{\JFM}[1]{JFM #1}

\theoremstyle{plain}
\newtheorem{theorem}{Theorem}
\newtheorem*{maintheorem}{Theorem~\ref{theorem:main}}
\newtheorem*{main2theorem}{Theorem~\ref{theorem:main2}}

\newtheorem{lemma}[theorem]{Lemma}
\newtheorem{proposition}[theorem]{Proposition}

\numberwithin{theorem}{section}
\numberwithin{conjecture}{section}
\numberwithin{problem}{section}

\numberwithin{equation}{section}

\newcommand{\R}{{\mathbb R}}
\newcommand{\Z}{{\mathbb Z}}

\newcommand{\Co}{\textup{Co}}
\newcommand{\ev}{\mathop{\textup{ev}}\nolimits}
\newcommand{\vol}{\mathop{\textup{vol}}\nolimits}
\newcommand{\Tr}{\mathop{\textup{Tr}}\nolimits}
\newcommand{\co}{\colon}

\title[Optimality and uniqueness of the Leech lattice]{Optimality and uniqueness of the Leech lattice among
lattices}

\dedicatory{Dedicated to Oded Schramm (10 December 1961 -- 1 September 2008)}

\author{Henry Cohn}
\address{Microsoft Research\\
One Microsoft Way\\
Redmond, WA 98052-6399, United States}
\curraddr{Microsoft Research New England\\
One Memorial Drive\\
Cambridge, MA 02142, United States}
\email{cohn@microsoft.com}

\author{Abhinav Kumar}
\address{Department of Mathematics\\
Harvard University\\
Cambridge, MA 02138, United States}
\curraddr{Department of Mathematics, Room 2-169\\
Massachusetts Institute of Technology\\
Cambridge, MA 02139, United States}
\email{abhinav@math.mit.edu}

\thanks{Kumar was supported by a summer internship
in the Theory Group at Microsoft Research and by a Putnam graduate
fellowship at Harvard University.\\
\ \\
Henry Cohn and Abhinav Kumar, \emph{Optimality and uniqueness of the Leech lattice among
lattices}, Annals of Mathematics \textbf{170} (2009), 1003--1050.}

\begin{document}

\begin{abstract}
We prove that the Leech lattice is the unique densest lattice in
$\R^{24}$.  The proof combines human reasoning with computer
verification of the properties of certain explicit polynomials. We
furthermore prove that no sphere packing in $\R^{24}$ can exceed the
Leech lattice's density by a factor of more than $1+1.65\cdot
10^{-30}$, and we give a new proof that $E_8$ is the unique densest
lattice in $\R^8$.
\end{abstract}

\maketitle

\section{Introduction}

It is a long-standing open problem in geometry and number theory
to find the densest lattice in $\R^n$.  Recall that a lattice
$\Lambda \subset \R^n$ is a discrete subgroup of rank $n$; a
minimal vector in $\Lambda$ is a nonzero vector of minimal
length.  Let $|\Lambda| = \vol(\R^n/\Lambda)$ denote the covolume
of $\Lambda$, i.e., the volume of a fundamental parallelotope or
the absolute value of the determinant of a basis of $\Lambda$.  If
$r$ is the minimal vector length of $\Lambda$, then spheres of
radius $r/2$ centered at the points of $\Lambda$ do not overlap
except tangentially.  This construction yields a sphere packing
of density
$$
\frac{\pi^{n/2}}{(n/2)!}
\left(\frac{r}{2}\right)^n\frac{1}{|\Lambda|},
$$
since the volume of a unit ball in $\R^n$ is $\pi^{n/2}/(n/2)!$,
where for odd $n$ we define $(n/2)! = \Gamma(n/2+1)$. The densest
lattice in $\R^n$ is the lattice for which this quantity is
maximized.

There might be several distinct densest lattices in the same
dimension.  For example, the greatest density known in $\R^{25}$
is achieved by at least $23$ distinct lattices, although they are
not known to be optimal.  (See pages xix and~178 of \cite{CS} for
the details.)  We will speak of ``the densest lattice'' because it
sounds more natural.

The problem of finding the densest lattice is a special case of the sphere
packing problem, but there is no reason to believe that the densest sphere
packing should come from a lattice.  In particular, in $\R^{10}$ the densest
packing known is the Best packing $P_{10c}$, which is not a lattice packing
(see \cite[p.~140]{CS}).  It is conjectured that lattices are suboptimal in
all sufficiently high dimensions. However, many of the most interesting
packings in low dimensions are lattice packings, and lattices have strong
connections with other fields such as number theory. For example, the Hermite
constant $\gamma_n$ is defined to be the smallest constant such that for
every positive-definite quadratic form $Q(x_1,\dots,x_n)$ of determinant $D$,
there is a nonzero vector $(v_1,\dots,v_n) \in \Z^n$ such that
$Q(v_1,\dots,v_n) \le \gamma_n D^{1/n}$.  Finding the maximum density of a
lattice packing in $\R^n$ is equivalent to computing $\gamma_n$.

The densest lattice in $\R^n$ is known for $n \le 8$: the answers are the
root lattices $A_1$, $A_2$, $A_3$, $D_4$, $D_5$, $E_6$, $E_7$, and $E_8$. For
$n=3$ this is due to Gauss \cite{G}, for $4 \le n \le 5$ to Korkine and
Zolotareff \cite{KZ1}, \cite{KZ2}, and for $6 \le n \le 8$ to Blichfeldt
\cite{Bl}. However, before the present paper no further cases had been solved
since 1935. In 1946 Chaundy claimed to have dealt with $n=9$ and $n=10$, but
his paper \cite{Ch} implicitly assumes that a densest lattice in $\R^n$ must
contain one in $\R^{n-1}$ as a cross section.  That is known to be false (see
\cite{lam}), so the paper appears irreparably flawed.

In each of the solved cases, the optimal lattice is furthermore
known to be unique, up to scaling and isometries. This was proved
simultaneously with the optimality for $n \le 5$,  for $n=6$ it
was proved by Barnes \cite{Ba}, and for $6 \le n \le 8$ it was
proved by Vet\v{c}inkin \cite{V}.

In this paper we deal with $n=24$ (the theorem numbering is as it
will appear later in the paper):

\begin{maintheorem}
The Leech lattice is the unique densest lattice in $\R^{24}$, up
to scaling and isometries of $\R^{24}$.
\end{maintheorem}

In terms of the Hermite constant, $\gamma_{24}=4$.  We also give
a new proof for $E_8$:

\begin{main2theorem}[Blichfeldt, Vet\v{c}inkin]
The $E_8$ root lattice is the unique densest lattice in $\R^{8}$,
up to scaling and isometries of $\R^{8}$.
\end{main2theorem}

Our work is motivated by the paper \cite{CE} by Cohn and Elkies (see
also \cite{Co}), which proves upper bounds for the sphere packing
density.  In particular, the main theorem in \cite{CE} is an analogue
for sphere packing of the linear programming bounds for
error-correcting codes: given a function satisfying certain linear
inequalities one can deduce a density bound.  It is not known how to
choose the function to optimize the bound, but in $\R^8$ and $\R^{24}$
one can come exceedingly close to the densities of $E_8$ and the Leech
lattice, respectively.  This observation, together with analogies with
error-correcting codes and spherical codes, led Cohn and Elkies to
conjecture that their bound is sharp in $\R^8$ and $\R^{24}$, which
would solve the sphere packing problem in those dimensions.  While we
cannot yet fully carry out that program, in this paper we show how to
combine the methods of \cite{CE} with results on lattices and
combinatorics to deal with the special case of lattice packings. We
will deal primarily with the Leech lattice, because that case is new
and more difficult, but in Section~\ref{section:E8} we will discuss
$E_8$.

One might hope to use a relatively simple method. Section~8 of
\cite{CE} shows how to prove that the Leech lattice is the unique
densest periodic packing (i.e., union of finitely many translates
of a lattice) in $\R^{24}$, if a function from $\R^{24}$ to $\R$
with certain properties exists.  Using a computer, one can find
functions that very nearly have those properties, and the
techniques from Section~8 of \cite{CE} can then be used
straightforwardly to prove an approximate version of the
uniqueness result: every lattice that is at least as dense as the
Leech lattice must be close to it. It is known (see
\cite[p.~176]{M}) that the Leech lattice is a strict local
optimum for density, among lattices. Thus, if one can show that
every denser lattice is sufficiently close, then it proves that
the Leech lattice is the unique densest lattice in $\R^{24}$.

Unfortunately, this approach seems completely infeasible if
carried out in the most straightforward way. When one naively
imitates the techniques from Section~8 of \cite{CE} in an
approximate setting, one loses a tremendous factor in the bounds,
and that puts the required computer searches far beyond what we
are capable of.  In this paper we salvage the approach by using
more sophisticated arguments that take advantage of special
properties of the Leech lattice.  In particular, we make use of
three beautiful facts about the Leech lattice: its automorphism
group acts transitively on pairs of minimal vectors with the same
inner product, its minimal vectors form an association scheme
when pairs are grouped according to their inner products, and its
minimal vectors form a spherical $11$-design.  (Note that the
second property follows from the first, but our work uses another
proof of it, from \cite{DGS}.)

Our proof depends on computer calculations in some places, but they can
be carried out relatively quickly, in less than one hour using a
personal computer. The calculations are all done using exact arithmetic
and are thus rigorous.  We have fully documented all of our
calculations and made available commented code for use in checking the
results or carrying out further investigations.  See
Appendix~\ref{appendix:computer} for details.

Appendix~\ref{appendix:background} contains very brief introductions to
several topics: the Leech lattice, linear programming bounds, spherical
designs, and association schemes. For more details, see \cite{CS}.  The
expository articles \cite{Ea} and \cite{Eb} also provide useful background
and context, although they do not include everything we need.

\section{Outline of proof}

We wish to show that the Leech lattice, henceforth denoted by
$\Lambda_{24}$, is the unique densest sphere packing among all
lattices in $\R^{24}$.  Let $\Lambda$ be any lattice in $\R^{24}$
that is at least as dense as $\Lambda_{24}$.  Without loss of
generality, we assume that $\Lambda$ has covolume $1$. Then the
restriction on its density simply means its minimal vectors have
length at least $2$.

We first show, using linear programming bounds, that $\Lambda$
has exactly $196560$ vectors of length approximately $2$ (called
nearly minimal vectors), and that the next smallest vector length
is approximately $\sqrt{6}$.

We rescale those $196560$ nearly minimal vectors to lie on the unit
sphere. Then they form a spherical code with minimal angle at least
$\varphi$, where $\cos\varphi$ is very near to (and greater than or
equal to) $1/2$. Note that in $S^{23}$ there is a unique spherical code
of this size with minimal angle $\pi/3 = \cos^{-1}(1/2)$, and it is the
kissing arrangement of $\Lambda_{24}$; the spherical code derived from
$\Lambda$ should be a small perturbation of this configuration.

Using linear programming bounds, we show that the inner products
of the unit vectors are approximately $0, \pm 1/4, \pm 1/2, \pm
1$.  We prove that if pairs of vectors are grouped according to
their inner products, then they form an association scheme with
the same valencies and intersection numbers as in the case of
$\Lambda_{24}$, and that it must therefore be the same
association scheme.  This isomorphism gives us a correspondence
between minimal vectors of $\Lambda_{24}$ and nearly minimal
vectors of $\Lambda$, such that corresponding inner products are
approximately equal.

Using this correspondence, we find a basis of $\Lambda$ whose
Gram matrix is close to the Gram matrix of a basis of
$\Lambda_{24}$.  Finally, from the strict local optimality of the
Leech lattice we conclude that $\Lambda$ must in fact be the
Leech lattice.

\section{Notation}\label{section:notation}

We begin by recording our normalizations of some special functions (which are
always as in \cite{AAR}), and by defining some notation.

The Laguerre polynomials $L_i^\alpha(z)$ are defined by the
initial conditions $L_0^\alpha(z) = 1$ and $L_1^\alpha(z) =
1+\alpha-z$ and the recurrence
$$
iL_{i}^\alpha(z) =
(2i-1+\alpha-z)L_{i-1}^\alpha(z)-(i+\alpha-1)L_{i-2}^\alpha(z)
$$
for $i \ge 2$. They are orthogonal polynomials with respect to the
measure $e^{-x}x^\alpha \, dx$ on $[0,\infty)$. If $\alpha = n/2-1$,
then the functions on $\mathbb{R}^n$ given by $x \mapsto e^{-\pi |x|^2}
L_i^\alpha(2\pi |x|^2)$ form an orthogonal basis of the radial
functions in $L^2(\R^n)$, and they are also eigenfunctions of the
Fourier transform with eigenvalue $(-1)^i$ (see (4.20.3) in
\cite{Leb}).  Here we normalize the Fourier transform by
$$
\widehat{f}(t) = \int_{\R^n} f(x) e^{2\pi i \langle x,t \rangle}
\, dx.
$$
Note also that with this normalization of the Fourier transform,
the Poisson summation formula states that
$$
\sum_{x \in \Lambda} f(x) = \frac{1}{|\Lambda|} \sum_{t \in
\Lambda^*} \widehat{f}(t),
$$
if $f \co \R^n \to \R$ is a Schwartz function, $\Lambda \subset \R^n$
is a lattice, and
$$
\Lambda^* = \{ y \in \R^n : \langle x,y \rangle \in \Z \textup{
for all $x \in \Lambda$} \}
$$
is its dual lattice.  (See (28) in \cite[p.~49]{K}.)

The ultraspherical (or Gegenbauer) polynomials $C_i^\lambda(z)$
are defined by the initial conditions $C_0^\lambda(z)=1$ and
$C_1^\lambda(z) = 2\lambda z$ and the recurrence
$$
i C^\lambda_i(z) = 2(i+\lambda-1)z C^\lambda_{i-1}(z) -
(i+2\lambda-2)C^\lambda_{i-2}(z)
$$
for $i \ge 2$. They are orthogonal polynomials with respect to
the measure
$$
(1-x^2)^{\lambda-1/2}\,dx
$$
on $[-1,1]$.  When
$\lambda = n/2-1$, that measure is proportional to the projection
of the surface measure from $S^{n-1}$ onto an axis, and the
ultraspherical polynomials play a fundamental role in the theory
of spherical harmonics in $\R^n$.  Up to scaling, the
ultraspherical polynomial $C_i^\lambda$ is the same as the Jacobi
polynomial $P_i^{(\alpha,\alpha)}$, where $\alpha = \lambda-1/2$.

Whenever we use Laguerre or ultraspherical polynomials, we will
always set $\alpha=\lambda=n/2-1$, where $n=24$ in the Leech
lattice proof and $n=8$ in the $E_8$ proof.  The term
``ultraspherical coefficient'' will mean a coefficient in the
expansion of a polynomial as a linear combination of
ultraspherical polynomials.

Throughout this paper, $\Lambda_{24}$ will denote the Leech lattice, and
$\Lambda$ will denote any lattice in $\R^{24}$ that is at least as dense and
satisfies $|\Lambda|=1$ (except in Section~\ref{section:localopt}, where
$\Lambda$ denotes an arbitrary lattice, and in Section~\ref{section:E8},
which deals with $E_8$). We think of $\Lambda$ as being an optimal lattice,
but we will not use that assumption.  We do not even need to know \emph{a
priori} that a global optimum for density is achieved, although
\cite[\S17.5]{GL} shows that it is.

Whenever $f \co \R^n \to \R$ is a radial function and $r \in
[0,\infty)$, we will write $f(r)$ for the common value $f(x)$ with $x
\in \R^n$ satisfying $|x|=r$.

The surface volume of the unit sphere $S^{n-1} \subset \R^n$ will
be denoted by
$$
\vol(S^{n-1}) = n \frac{\pi^{n/2}}{(n/2)!}.
$$
It is important to keep in mind that $\vol(S^{n-1})$ is not the
volume of the enclosed ball.

We will make use of the numerical values
\begin{align*}
\varepsilon &=  6.733\cdot 10^{-27},\\
\mu &= 3.981 \cdot 10^{-13},\\
\nu &= 3.219 \cdot 10^{-12}, \textup{ and}\\
\omega &= 1.703 \cdot 10^{-11}
\end{align*}
throughout the paper.  Each will be defined the first time it is
used, but we have collected the values here for easy reference.
Note that terminating decimal expansions such as these represent
exact rational numbers, not floating point approximations.

\section{Nearly minimal vectors}\label{section:nearlyminimal}

In this section, we show that $\Lambda$ must have exactly
$196560$ vectors of length near $2$ (this will be made more
precise below).  The first subsection examines which vector
lengths are possible in $\Lambda$, and the second then counts the
nearly minimal vectors.

\subsection{Restrictions on the lengths of vectors}\label{howlarge}

Let $f \co \R^{24} \to \R$ be a radial Schwartz function with the
following properties: $f(0) = \widehat{f}(0)=1$, $f(x) \le 0$ for $|x|
\ge r$ (for some number $r$), and $\widehat{f}(t) \ge 0$ for all $t$.
Proposition~3.2 of \cite{CE} says that if such a function exists, then
the sphere packing density in $\R^{24}$ is bounded above by
$$
\frac{\pi^{12}}{12!}\left(\frac{r}{2}\right)^{24}.
$$
If we could find such a function with $r=2$, then it would prove
that $\Lambda_{24}$ has the greatest density among all sphere
packings in $\R^{24}$, not just lattice packings.  Cohn and Elkies
conjecture that such a function exists, but the best they achieve
in \cite{CE} is $r \le 2\cdot1.00002946$.

Our proof begins by constructing an explicit function $f$ with
$$
r \le 2\big(1+6.851\cdot10^{-32}\big).
$$
Note that the existence of such a function proves that no sphere
packing in $\R^{24}$ can have density greater than
$1+1.65\cdot10^{-30}$ times the density of $\Lambda_{24}$.

Unfortunately, the function we construct is extremely complicated
(it would take far too much space to write it down here). It
consists of a polynomial of degree $803$ with rational
coefficients, evaluated at $2\pi|x|^2$ and multiplied by
$e^{-\pi|x|^2}$. It was constructed by a lengthy computer
calculation to optimize the value of $r$ using Newton's method,
and even verifying that it has the properties used below requires
a computer, although fortunately that is much easier than finding
the function. The accompanying computer file \texttt{verifyf.txt}
includes code to verify all the assertions about $f$ in this
subsection of the paper.  See Appendix~\ref{appendix:computer}
for details.

We can use techniques similar to those in \cite{CE} to study the
size of the short vectors in $\Lambda$.  First, we need two
lemmas. Set
$$
\varepsilon = 6.733\cdot 10^{-27},
$$
and call a nonzero vector in $\Lambda$ \textit{nearly minimal\/} if it
has length at most $2(1+\varepsilon)$.  The reason for this choice of
$\varepsilon$ will be apparent from Proposition~\ref{lengthbounds}
below.

Consider what happens if we rescale the nearly minimal vectors so that
they all lie on the unit sphere. These vectors determine a spherical
code
$$
\mathcal{C}_\Lambda = \{ u/|u| : u \textrm{ a nearly minimal
vector} \}
$$
on the unit sphere $S^{23}$, and the following lemma bounds its
minimal angle.  (See Appendix~\ref{appendix:background} for
background on spherical codes.)

\begin{lemma}\label{lemma:angular}
If $u$ and $v$ are nearly minimal vectors with $u \ne v$, then
the angle $\varphi$ between $u$ and $v$ satisfies
$$
\cos\varphi \le 1 - \frac{1}{2(1+\varepsilon)^2}.
$$
\end{lemma}

\begin{proof}
We have $|u|, |v| \in [2,2(1+\varepsilon)]$ and $|u-v| \ge 2$. By
the law of cosines,
$$
\cos\varphi = \frac{|u|^2+|v|^2-|u-v|^2}{2|u||v|} \le
\frac{|u|^2+|v|^2-4}{2|u||v|}.
$$
The bound $(|u|^2+|v|^2-4)/(2|u||v|)$ is convex as a function of $|u|$
and $|v|$ individually, and hence it is maximized at one of the
vertices of the square $[2,2(1+\varepsilon)]^2$.  In fact, the maximum
occurs when $|u|=|v| = 2(1+\varepsilon)$, in which case the bound
becomes
$$
1 - \frac{1}{2(1+\varepsilon)^2}.
$$
\end{proof}

\begin{lemma}\label{atmost196560}
There are at most $196560$ nearly minimal vectors in $\Lambda$.
\end{lemma}

\begin{proof}
This lemma is a straightforward application of the linear programming
bounds for spherical codes (see Chapter~9 of \cite{CS}, or
Appendix~\ref{appendix:background} for a brief summary). Let
$$
f_\varepsilon(x) = K_\varepsilon
(x+1)\left(x+\frac{1}{2}\right)^2\left(x+\frac{1}{4}\right)^2
x^2\left(x-\frac{1}{4}\right)^2
\left(x-\left(1-{\frac{1}{2(1+\varepsilon)^2}}\right)\right),
$$
where the constant $K_\varepsilon$ is chosen so that $f_\varepsilon$
has zeroth ultraspherical coefficient $1$.  (The normalization is
irrelevant for this proof, but it will be important later in the paper,
so we use it here for consistency.) If $\varepsilon$ were 0, this
polynomial would be the one used to solve the kissing problem exactly
in $\R^{24}$ (see Chapter~13 of \cite{CS}). With the current value of
$\varepsilon$, the polynomial $f_\varepsilon$ has nonnegative
ultraspherical coefficients and still proves that there are fewer than
$196561$ spheres in any spherical code in $\R^{24}$ with minimal angle
as in Lemma~\ref{lemma:angular}.  (We check this assertion in the
computer file \texttt{verifyrest.txt}.  In fact, the bound is less than
$196560+10^{-19}$.) Thus, there can be at most $196560$ nearly minimal
vectors in $\Lambda$.
\end{proof}

In addition to the definition of $\varepsilon$ above, set
\begin{align*}
\mu &= 3.981 \cdot 10^{-13},\\
\nu &= 3.219 \cdot 10^{-12}, \textup{ and}\\
\omega &= 1.703 \cdot 10^{-11}.
\end{align*}

\begin{proposition} \label{lengthbounds}
Every nonzero vector in $\Lambda$ has length in
$$
\big[2,2(1+\varepsilon)\big) \cup
\big(\sqrt{6}(1-\mu),\sqrt{6}(1+\mu)\big) \cup
\big(\sqrt{8}(1-\nu),\sqrt{8}(1+\nu)\big) \cup
\big(\sqrt{10}(1-\omega),\infty\big).
$$
\end{proposition}

\begin{proof}
By the Poisson summation formula,
$$
\sum_{x \in \Lambda} f(x) = \sum_{t\in\Lambda^*} \widehat{f}(t).
$$
Because $\Lambda$ is at least as dense as $\Lambda_{24}$ (and has
covolume~$1$), all nonzero vectors $x \in \Lambda$ satisfy $|x|
\ge 2$.  Because $r \le 2(1+\varepsilon)$, by
Lemma~\ref{atmost196560} there can be at most $196560$ vectors in
$\Lambda$ with $2 \le |x| \le r$. Within that range, $f$ is a
decreasing function of the radius, and we have
$$
196560 f(2) < 1.644104221 \cdot 10^{-30}
$$
(recall from Section~\ref{section:notation} that $f(2)$ denotes
the common value $f(x)$ with $|x|=2$).

The key properties of $f$ are $f(0) = \widehat{f}(0)=1$, $f(x) \le 0$
for $|x| \ge r$, and $\widehat{f}(t) \ge 0$ for all $t$. It follows
that
$$
1 + 1.644104221 \cdot 10^{-30} + \sum_{x \in \mathcal{N}} f(x) \ge 1,
$$
where $\mathcal{N}$ is the set of vectors in $\Lambda$ at which
$f$ is negative.  No vector in $\Lambda$ can occur within any
region on which $f$ is less than $(1/2)(-1.644104221) \cdot
10^{-30}$; the extra factor of $2$ comes from the fact that $f(-x)
= f(x)$ since $f$ is a radial function (if $x \in \mathcal{N}$
then $-x \in \mathcal{N}$ as well).  That rules out all radii in
the set
$$
\big[2(1+\varepsilon), \sqrt{6}(1-\mu)\big] \cup
\big[\sqrt{6}(1+\mu),\sqrt{8}(1-\nu)\big]\cup
\big[\sqrt{8}(1+\nu),\sqrt{10}(1-\omega)\big].
$$
In the computer file \texttt{verifyf.txt} we prove this by
examining the radial derivative of $f$.
\end{proof}

Much of the rest of the proof would still work if $\varepsilon$,
$\mu$, $\nu$, and $\omega$ were somewhat larger.  The two main
places where they must be small are the final inequality
\eqref{eq:final} and the intersection number calculations in
Subsection~\ref{subsec:intnum} (as well as the bounds used there).
In each case they could be slightly larger, but not by a factor
of $100$.

\subsection{196560 nearly minimal vectors.}

We can now show that there are exactly $196560$ nearly minimal
vectors.  We know from Lemma~\ref{atmost196560} that there are at
most $196560$ of them.  For the other direction, a lower bound
greater than $196559$, we need a new kind of linear programming
bound.  Recall that we have shown that all nonzero vectors in
$\Lambda$ are either nearly minimal or have lengths greater than
$\sqrt{6}(1-\mu)$.

Suppose we knew that all nonzero vectors either have length exactly $2$
or at least $\sqrt{6}$.  (We will first explain our method under these
overly optimistic hypotheses. Lemma~\ref{atleast196559} will then apply
it using the actual bounds we have proved.) One might hope to count the
nearly minimal vectors using a Schwartz function $g \co \R^{24} \to \R$
such that $g(x) \le 0$ for $|x| \ge \sqrt{6}$, $\widehat{g}(t) \ge 0$
for all $t$, and $g(2)
> 0$. Given such a function, Poisson summation implies that
$$
\sum_{x \in \Lambda} g(x) = \sum_{t \in \Lambda^*} \widehat{g}(t),
$$
and hence
$$
g(0) + N g(2) \ge \widehat{g}(0)
$$
if there are $N$ minimal vectors.  Thus, $N \ge
(\widehat{g}(0)-g(0))/g(2)$.  We conjecture that $g$ can be
chosen so that $(\widehat{g}(0)-g(0))/g(2) = 196560$, which is
the largest possible value because the Leech lattice has $196560$
minimal vectors. We can construct functions that come quite close
to this bound, and will use one of them to prove the following
lemma.

\begin{lemma}\label{atleast196559}
There are more than $196559$ nearly minimal vectors in $\Lambda$.
\end{lemma}

\begin{proof}
Define $z_1,\dots,z_{10}$ by $z_i = \lfloor 4\pi(i+1)10^8\rfloor/10^8$.
In other words, they have the following values:

\begin{center}
\begin{tabular}{c|ccccc}
$i$ & 1 & 2 & 3 & 4 & 5\\
\hline $z_i$ & 25.13274122 & 37.69911184 & 50.26548245 &
62.83185307 & 75.39822368
\end{tabular}

\begin{tabular}{c|ccccc}
$i$ & 6 & 7 & 8 & 9 & 10\\
\hline $z_i$ & 87.96459430 & 100.53096491 & 113.09733552 &
125.66370614& 138.23007675
\end{tabular}
\end{center}

There are unique coefficients
$a_1,\dots,a_{37}$ such that
$$
1+\sum_{i=1}^{37} a_i L_i^{11}(x)
$$
has a single root at $z_2$ and a double root at $z_i$ for $i \geq
3$, and
$$
1+\sum_{i=1}^{37} (-1)^i a_i L_i^{11}(x)
$$
has a double root at each $z_i$ for $i \geq 1$.  Neither
polynomial has any other nonnegative roots.  (We check this in
the computer file \texttt{verifyg.txt} using Sturm's theorem,
except that we do not check the uniqueness of the coefficients
because we do not require it.)

Define $g \co \R^{24} \to \R$ by
$$
g(x) = \left(1+\sum_{i=1}^{37} a_i
L_i^{11}(2\pi|x|^2)\right)e^{-\pi|x|^2}.
$$
It follows that
$$
\widehat{g}(t) = \left(1+\sum_{i=1}^{37} (-1)^i a_i
L_i^{11}(2\pi|x|^2)\right)e^{-\pi|x|^2}.
$$
We have $g(x) \le 0$ for $|x| \ge \sqrt{6}(1-\mu)$, because $z_2
< 2\pi \cdot 6(1-\mu)^2$, the function $g$ changes sign only at
$z_2$, and $g(0)>0$. For all $t$, we have $\widehat{g}(t) \ge 0$.

Apply Poisson summation to $g$, to deduce
$$
\sum_{x \in \Lambda} g(x) = \sum_{t \in \Lambda^*} \widehat{g}(t).
$$
Applying the two inequalities above shows that
$$
g(0) + \sum_{x \in \mathcal{M}} g(x) \ge \widehat{g}(0),
$$
where $\mathcal{M}$ is the set of nearly minimal vectors.  The
function $g(x)$ is positive and a decreasing function of $|x|$ on
the interval $[2,2(1+\varepsilon)]$, so
$$
g(0) + |\mathcal{M}| g(2) \ge \widehat{g}(0).
$$
However,
$$
\frac{\widehat{g}(0) - g(0)}{g(2)} > 196559,
$$
so there are more than $196559$ nearly minimal vectors. (All
these inequalities are checked in \texttt{verifyg.txt}.)
\end{proof}

Thus, there must be exactly $196560$ nearly minimal vectors, as
desired.  We conjecture that this method could be used to recover
each of the coefficients of the Leech lattice's theta series, but
we will not need that for our proof.

\section{Inner products in the spherical code}\label{section:sphericalcode}

We now continue to study the polynomial
$$
f_\varepsilon(x) =
K_\varepsilon(x+1)\left(x+\frac{1}{2}\right)^2\left(x+\frac{1}{4}\right)^2
x^2\left(x-\frac{1}{4}\right)^2
\left(x-\left(1-{\frac{1}{2(1+\varepsilon)^2}}\right)\right)
$$
from the previous section.

Note that
$$
1-\frac{1}{2(1+\varepsilon)^2} < \frac{1}{2} + \varepsilon.
$$
Thus, all the inner products except $1$ in the spherical code
$\mathcal{C}_\Lambda$ are at most $1/2 + \varepsilon$.  We seek bounds for
how far from $0, \pm 1/4, \pm 1/2, \pm 1$ they can be.  The $\pm 1$ cases
must be exact, because of Lemma~\ref{lemma:angular} and the fact that $u \in
\mathcal{C}_\Lambda$ if and only if $-u \in \mathcal{C}_\Lambda$ (i.e., the
code is antipodal).

Because the zeroth ultraspherical coefficient of $f_\varepsilon$ is
$1$, it follows from the usual proof of the linear programming bounds
for spherical codes (see Appendix~\ref{appendix:background}) that
$$
196560f_{\varepsilon}(1) + \sum_{x\ne y} f_\varepsilon(\langle x,
y \rangle) \ge 196560^2,
$$
where the sum is over vectors in the spherical code.  Because the
four inner products $\langle x,y \rangle$, $\langle -x,-y
\rangle$, $\langle y,x\rangle$, and $\langle -y,-x \rangle$ are
all equal, and all the terms in the sum are nonpositive, we see
that no term in the sum can be less than
$$
\frac{196560^2-196560f_\varepsilon(1)}{4},
$$
which is approximately $-1.99\cdot 10^{-15}$.

Now a short calculation implies that all the inner products must
be within $6.411 \cdot 10^{-9}$ of one of the numbers $0, \pm
1/4, \pm 1/2, \pm 1$.  The exponent is only $9$ because
$f_\varepsilon$ has double roots, and one can stray quite far
from a double root without substantially changing the function's
value.

In the rest of the paper, we make the following definition: let
$$
\sigma = \max_{u,v \in \mathcal{C}_\Lambda} \min \{|\langle u,v
\rangle - \alpha| : \alpha \in S\},
$$
where $S = \{0, \pm 1/4, \pm 1/2, \pm 1\}$.  In other words,
$\sigma$ is the maximum ``error'' in the inner products. We have
just shown that $\sigma \leq 6.411 \cdot 10^{-9}$.  We will
improve the upper bound for $\sigma$ substantially, and ultimately
we will show $\sigma = 0$.

\subsection{Better bounds for $\sigma$}

Recall that we computed in Section~\ref{howlarge} that vectors of
$\Lambda$ with length close to $\sqrt{6}$ must have length in the
interval $\sqrt{6}(1-\mu), \sqrt{6}(1+\mu)$. Therefore if we have
nearly minimal vectors $u, v$  with $\langle u, v \rangle \approx
1$ (with error less than $10^{-3}$, say), then we see that $|u-v|
\approx \sqrt{6}$. Therefore
$$
6(1-\mu)^2 \leq \langle u, u \rangle + \langle v,v \rangle -2
\langle u,v \rangle  \leq 6(1+\mu)^2.
$$
In addition we have
$$
4 \leq \langle u,u \rangle \leq 4(1+\varepsilon)^2.
$$
Therefore
$$
8 - 6(1+\mu)^2 \leq  2\langle u,v \rangle \leq 8(1+\varepsilon)^2
-6(1-\mu)^2,
$$
so
$$
\frac{4-3(1+\mu)^2}{4(1+\varepsilon)^2} \leq \left\langle \frac{u}{|u|},
\frac{v}{|v|} \right\rangle \leq (1+\varepsilon)^2 - \frac{3}{4}(1-\mu)^2.
$$
Similarly we get for $\langle u, v \rangle \approx 0$ that
$$
\frac{1-(1+\nu)^2}{(1+\varepsilon)^2} \leq \left\langle \frac{u}{|u|},
\frac{v}{|v|} \right\rangle \leq (1+\varepsilon)^2 - (1-\nu)^2,
$$
and for $\langle u,v \rangle \approx 2$ that
$$
\frac{2-(1+\varepsilon)^2}{2(1+\varepsilon)^2} \leq \left\langle \frac{u}{|u|},
\frac{v}{|v|} \right\rangle \leq (1+\varepsilon)^2 - \frac{1}{2}
$$
because $2 \leq |u-v| \leq 2(1+\varepsilon)$. Combining these (and if
necessary replacing $v$ with $-v$), we find that for $u \neq \pm v$,
the inner product $\langle u/|u|, v/|v| \rangle$ differs from an
element of $\{0, \pm 1/4, \pm 1/2\}$ by at most $6.43801 \cdot
10^{-12}$.  We conclude that the error $\sigma$ in the inner products
is at most $6.43801 \cdot 10^{-12}$.

However, we will be able to get much better bounds in
Section~\ref{section:innerbounds}, once we have shown that the
spherical code $\mathcal{C}_\Lambda$ gives us an association scheme.

\section{Association schemes}\label{section:associationschemes}

We would like to turn the $196560$ points in the spherical code
$\mathcal{C}_\Lambda$ into a $6$-class association scheme by grouping
pairs according to their approximate inner products. (See
Appendix~\ref{appendix:background} for background on association
schemes.)  It is not clear that this in fact defines an association
scheme, but we will show that it does. Furthermore, we will show that
this association scheme is isomorphic to the one derived from
$\Lambda_{24}$. To achieve this, we show that the intersection numbers
are the same as in $\Lambda_{24}$. That will also show that it is an
association scheme, by showing that the intersection numbers are
independent of the pair of points. We use the same techniques as
\cite{DGS}, but we need to keep track of error bounds.

\subsection{Spherical design}

First, we show that $\mathcal{C}_\Lambda$ is nearly a spherical
$10$-design. (See Appendix~\ref{appendix:background} for
background on spherical designs.) Let
$$
C_i(x) = \frac{C_i^{11}(x)}{C_i^{11}(1)} \cdot
\frac{\binom{22+i}{22}+\binom{21+i}{22}}{\vol(S^{23})}.
$$
The advantage of this normalization of the ultraspherical
polynomials is that for every finite subset $\mathcal{C}$ of
$S^{23}$,
$$
\sum_{x,y \in \mathcal{C}} C_i(\langle x,y\rangle) = \left|
\sum_{z \in \mathcal{C}} \ev_i(z) \right|^2,
$$
where $\ev_i(z)$ denotes the evaluation at $z$ map in the dual
space to the $i$-th degree spherical harmonics.  Although this
fact is well known (e.g., it is equivalent to Theorem~9.6.3 in
\cite{AAR}), we will explain it here for completeness, because
the correct normalization is important for our application.

\begin{lemma}\label{lemma:nonneg}
If $\mathcal{C}$ is a finite subset of $S^{23}$, then
$$
\sum_{x,y \in \mathcal{C}} C_i(\langle x,y\rangle) = \left|
\sum_{z \in \mathcal{C}} \ev_i(z) \right|^2.
$$
\end{lemma}

\begin{proof}
Let $d$ be the dimension of the space of spherical harmonics of
degree $i$, and let $S_1,\dots,S_d$ be an orthonormal basis of
that space.  By Theorem~9.6.3 of \cite{AAR},
$$
\sum_{j=1}^d S_j(w) S_j(z) = C_i(\langle w,z \rangle).
$$
Let $f = \sum_j a_j S_j$ be any spherical harmonic of degree
$i$.  Then
\begin{align*}
(\ev_i(z))(f) &= \sum_{j=1}^d a_j S_j(z)\\
&= \sum_{j=1}^d S_j(z)\int_{S^{23}} S_j(w) f(w) \, dw\\
&= \int_{S^{23}} \left(\sum_{j=1}^d S_j(w) S_j(z)\right) f(w) \,
dw\\
&= \int_{S^{23}} C_i(\langle w,z \rangle) f(w) \, dw.
\end{align*}
Thus,
$$
\left( \sum_{z \in \mathcal{C}} \ev_i(z) \right)(f) =
\int_{S^{23}} \left(\sum_{z \in \mathcal{C}} C_i(\langle w,z
\rangle)\right) f(w) \, dw.
$$
In other words, applying the element $\sum_{z \in \mathcal{C}}
\ev_i(z)$ of the dual space is the same as taking the inner
product with
$$
w \mapsto \sum_{z \in \mathcal{C}} C_i(\langle w,z \rangle).
$$
It follows that
\begin{align*}
\left| \sum_{z \in \mathcal{C}} \ev_i(z) \right|^2 &=
\int_{S^{23}} \left(\sum_{z \in \mathcal{C}} C_i(\langle w,z
\rangle)\right)^2 \, dw\\
&= \sum_{x,y \in \mathcal{C}} \int_{S^{23}} C_i(\langle w,x
\rangle) C_i(\langle w,y \rangle) \, dw\\
&= \sum_{x,y \in \mathcal{C}} (\ev_i(x)) (w \mapsto C_i(\langle
w,y \rangle))\\
&= \sum_{x,y \in \mathcal{C}} C_i(\langle x,y \rangle),
\end{align*}
as desired.
\end{proof}

The following lemma asserts that $\mathcal{C}_\Lambda$ is nearly
a spherical $10$-design:

\begin{lemma}\label{nearlydesign}
If $g \co S^{23} \to \R$ is a polynomial of total degree at most $10$,
then
$$
\left| \sum_{z \in \mathcal{C}_\Lambda} g(z) -
\frac{196560}{\vol(S^{23})} \int_{S^{23}} g(z) \, dz \right| \le
2.50193 \cdot 10^{-5} |g|_2,
$$
where $|g|_2$ denotes the norm on $L^2(S^{23})$.
\end{lemma}

\begin{proof}
Without loss of generality we may assume that $g$ is a harmonic
polynomial (for every polynomial on $\R^{24}$, there is a harmonic
polynomial of equal or lesser degree with the same restriction to
$S^{23}$; see equation (5) in \cite[p.~17]{H}). We will use the
polynomial
$$
f_\varepsilon(x) =
K_\varepsilon (x+1)\left(x+\frac{1}{2}\right)^2\left(x+\frac{1}{4}\right)^2
x^2\left(x-\frac{1}{4}\right)^2
\left(x-\left(1-{\frac{1}{2(1+\varepsilon)^2}}\right)\right)
$$
from earlier in the paper.  Recall that $f_\varepsilon$ is normalized
to have zeroth ultraspherical coefficient $1$  (using the standard
normalization of the ultraspherical polynomials, not the new
normalization $C_i$).  If $c_i$ denotes the coefficient of $C_i$ in
$f_\varepsilon$, then
\begin{align*}
196560f_\varepsilon(1) &\ge \sum_{x,y \in \mathcal{C}_\Lambda}
f_\varepsilon(\langle x,y \rangle)\\
&=  196560^2 + \sum_{i=1}^{10} c_i \sum_{x,y \in \mathcal{C}} C_i(\langle x,y\rangle)\\
&=  196560^2 + \sum_{i=1}^{10} c_i \left|\sum_{z \in
\mathcal{C}_\Lambda} \ev_i(z)\right|^2,
\end{align*}
from which it follows that
$$
\sum_{i=1}^{10} c_i \left|\sum_{z \in \mathcal{C}_\Lambda}
\ev_i(z)\right|^2 \le 196560f_\varepsilon(1)-196560^2 < 7.9775
\cdot 10^{-15}.
$$
Therefore,
$$
\sum_{i=1}^{10}\left| \sum_{z \in \mathcal{C}_\Lambda} \ev_i(z)
\right|^2 \le \left(7.9775 \cdot 10^{-15}\right) \max_i
\frac{1}{c_i} < 6.25964 \cdot 10^{-10}.
$$

We can now bound
$$
\left| \sum_{z \in \mathcal{C}_\Lambda} g(z) -
\frac{196560}{\vol(S^{23})} \int_{S^{23}} g(z) \, dz \right|.
$$
Write $g = \sum_{i=0}^{10} g_i$, where $g_i$ is homogeneous of
degree $i$.  The integral cancels with the $g_0$ term in the sum,
so we simply need to bound
$$
\left|\sum_{i=1}^{10} \sum_{z \in \mathcal{C}_\Lambda}
g_i(z)\right|.
$$
For that, we use the definition of the norm and the Cauchy-Schwarz
inequality to deduce that
\begin{align*}
\left|\sum_{i=1}^{10} \sum_{z \in \mathcal{C}_\Lambda}
g_i(z)\right| &\le \sum_{i=1}^{10} \left| \sum_{z \in
  \mathcal{C}_\Lambda} \ev_i(z) \right| |g_i|_2  \\
&\le \sqrt{\sum_{i=1}^{10}\left| \sum_{z \in
\mathcal{C}_\Lambda} \ev_i(z) \right|^2} \sqrt{\sum_{i=1}^{10}
|g_i|_2^2} \\\
 &\le 2.50193 \cdot 10^{-5} |g|_2,
\end{align*}
as desired.
\end{proof}

In fact, it follows immediately that $\mathcal{C}_\Lambda$ is
nearly a spherical $11$-design in the same sense, because it is
antipodal (if $x \in \mathcal{C}_\Lambda$ then $-x \in
\mathcal{C}_\Lambda$) and thus every homogeneous polynomial of
odd degree averages to $0$ over $\mathcal{C}_\Lambda$.  However,
we will not need that fact.

It is worth pointing out for completeness that the minimal vectors
in $\Lambda_{24}$ do not form a spherical $12$-design: if $y$ is a
minimal vector then the polynomial
$$
\left(16-\langle x,y \rangle^2\right)\left(4-\langle x,y \rangle^2
\right)^2 \left(1-\langle x,y \rangle^2\right)^2 \langle x,y
\rangle^2
$$
in $x$ vanishes at each minimal vector but does not average to
$0$ over the sphere of radius $2$, because it is
nonnegative on that sphere.

The constant $2.50193 \cdot 10^{-5}$ in Lemma~\ref{nearlydesign}
can be made somewhat smaller for polynomials of total degree at
most $8$ (using the same proof), and that is the only case we need
later.  However, the present bound suffices.

\subsection{Intersection numbers} \label{subsec:intnum}

We will now use the fact that $\mathcal{C}_\Lambda$ is nearly a
spherical $10$-design to determine the intersection numbers
(still following the techniques of \cite{DGS}). As a computational
aid, it is useful to know the following formula for averaging
homogeneous polynomials over the sphere:
$$
\frac{1}{\vol(S^{23})}\int_{S^{23}} g(z) \, dz =
\frac{\int_{\R^{24}} g(z) e^{-|z|^2} \, dz}{\int_0^\infty r^{\deg
g} e^{-r^2} \vol(S^{23}) r^{23} \, dr}.
$$

In the exact case of $\Lambda_{24}$, we can compute the intersection numbers
as follows.  For each $x,y \in \mathcal{C}_{24}$ (the spherical code derived
from the minimal vectors) with a specified inner product, we need to
determine the number of $z \in \mathcal{C}_{24}$ with specified inner
products with $x$ and $y$. Let $P_\gamma(\alpha,\beta)$ denote this number
when $\langle x,y \rangle = \gamma$, $\langle x,z \rangle = \alpha$, and
$\langle y,z \rangle = \beta$.  The cases $\gamma = \pm 1$ simply amount to
the valencies, which are determined automatically once the other intersection
numbers are determined. For instance $P_1(\alpha, \beta) = 0$ unless $\alpha
= \beta$, and $P_1(\alpha, \alpha) = \sum_{\beta} P_0(\alpha, \beta)$ once we
demonstrate that every vector in $\mathcal{C}_{24}$ has a vector in
$\mathcal{C}_{24}$ orthogonal to it (which follows from, say, showing that
$P_{\alpha}(0,0) \neq 0$ for each $\alpha \ne \pm 1$). Hence we will focus on
the remaining cases, i.e., $\gamma \neq \pm 1$.

For such a pair $(x, y)$ with $\langle x,y \rangle = \gamma$, \emph{a priori}
there are 49 unknowns $P_{\gamma}(\alpha, \beta)$ for $\alpha, \beta \in
\{-1,-1/2,-1/4,0,1/4,1/2,1\}$. We first note that $P_{\gamma}(1,\alpha) =
P_{\gamma}(\alpha,1) = \delta_{\alpha,\gamma}$, where $\delta$ is the
Kronecker delta. Similarly $P_{\gamma}(-1, \alpha) = P_{\gamma}(\alpha,-1) =
\delta_{\gamma, -\alpha}$. Thus we can eliminate $\pm 1$ from consideration,
which reduces the problem to finding only $25$ unknowns. We will find $25$
linear equations that determine these values.

Consider the polynomials $g_{i,j}(z) = \langle z, x \rangle^i
\langle z,y\rangle^j$ for $i,j  \in \{0,1,\dots,4\}$. Let $S$ be
the set $\{-1,-1/2,-1/4,0,1/4,1/2,1\}$. We then know that
$$
\sum_{\alpha,\beta \in S} \alpha^i \beta^j P_\gamma(\alpha,\beta)
= \frac{196560}{\vol(S^{23})} \int_{S^{23}} g_{i,j}(z) \, dz,
$$
because $\mathcal{C}_{24}$ is a spherical $11$-design (although even an
$8$-design would suffice). These equations can be solved to yield the unknown
values of $P_{\gamma}(\alpha,\beta)$. Note that the right-hand side does not
depend on the choice of $x$ and $y$, only on $\gamma= \langle x,y \rangle$.
Therefore we see that the solutions of these equations, which are the
intersection numbers, are independent of $x,y$ and only depend on $\gamma$.
We solve one such system of equations for each value of $\gamma$, and the
values of the intersection numbers are tabulated in
Table~\ref{table:intnums}, which is a complete list modulo the symmetries
$$
P_\gamma(\alpha,\beta) = P_\gamma(\beta,\alpha) =
P_\gamma(-\alpha,-\beta) = P_{-\gamma}(\alpha,-\beta).
$$
(These numbers are of course known, but they are not tabulated in
standard references such as \cite{CS}, so we record them here for
convenience.)

\begin{table}
\begin{center}
\begin{tabular}{lll}
$P_0(0,0) = 43164$ & $P_0(0,1/2) = 2464$ &
$P_0(0,1/4) = 22528$\\
$P_0(1/2,1/2) = 44$ & $P_0(1/2,1/4) = 1024$ & $P_0(1/4,1/4) =
11264$\\
\\
$P_{1/2}(0,0) = 49896$ & $P_{1/2}(0,1/2) = 891$  &
$P_{1/2}(0,1/4) = 20736$\\
$P_{1/2}(1/2,1/2) = 891$ & $P_{1/2}(1/2,-1/2) = 1$ &
$P_{1/2}(1/2,1/4) = 2816$\\
$P_{1/2}(1/2,-1/4) = 0$ & $P_{1/2}(1/4,1/4) = 20736$ &
$P_{1/2}(1/4,-1/4) = 2816$\\
\\
$P_{1/4}(0,0) = 44550$ & $P_{1/4}(0,1/2) = 2025$ &
$P_{1/4}(0,1/4) = 22275$\\
$P_{1/4}(1/2,1/2) = 275$ & $P_{1/4}(1/2,-1/2) = 0$ &
$P_{1/4}(1/2,1/4) = 2025$\\
$P_{1/4}(1/2,-1/4) = 275$ & $P_{1/4}(1/4,1/4) = 15400$ &
$P_{1/4}(1/4,-1/4) = 7128$\\
\\ 
\end{tabular}
\end{center}
\caption{Intersection numbers for the Leech lattice minimal
vectors.} \label{table:intnums}
\end{table}

What happens when we are not necessarily dealing with exactly
$\Lambda_{24}$, but rather with $\Lambda$? Suppose we have $x,y
\in \mathcal{C}_\Lambda$ and we want to determine the number of
$z \in \mathcal{C}_\Lambda$ with specified approximate inner
products with them.  Let $\widetilde{P}_\gamma(\alpha,\beta)$
denote the intersection numbers for $\mathcal{C}_\Lambda$ (which
may depend on $x$ and $y$);  here we use $\alpha$, $\beta$, and
$\gamma$ to denote exact elements of $\{0,\pm 1/4, \pm 1/2, \pm
1\}$, and the inner products from $\mathcal{C}_\Lambda$ are
required to be approximately equal to them. Then we have
\begin{align*}
\sum_{\alpha,\beta \in S} \alpha^i \beta^j
\widetilde{P}_{\gamma}(\alpha,\beta) &= \sum_{w \in
  \mathcal{C}_\Lambda, \langle w,x \rangle \approx \alpha, \langle w,y
  \rangle \approx \beta} \alpha^i \beta^j \\
&\approx  \sum_{w \in \mathcal{C}_\Lambda} \langle w,x \rangle^i
\langle w,y \rangle^j \\
&= \sum_{w \in
\mathcal{C}_\Lambda} g_{i,j}(w)\\
&\approx \frac{196560}{\vol(S^{23})} \int_{S^{23}}
g_{i,j}(z) \, dz\\
&\approx \frac{196560}{\vol(S^{23})} G_{i,j}(\gamma),
\end{align*}
where
$$
G_{i,j}(\gamma) = \int_{S^{23}} \langle z, u \rangle^i \langle
z,v\rangle^j\, dz
$$
with $\langle u,v \rangle = \gamma$ (recall that $\langle x,y
\rangle \approx \gamma$).

If these approximations are close enough, we should get about the
same values for $\widetilde{P}_{\gamma}(\alpha, \beta)$ as we did
for $P_{\gamma}(\alpha, \beta)$, and this will show that the
intersection numbers must be the same.

\begin{lemma}\label{lemma:approx}
Let $\alpha, \beta \in \{0,\pm 1/4, \pm 1/2, \pm 1\}$ and $a, b \in
[-(1/2+\sigma),1/2+\sigma] \cup \{\pm 1\}$ with $\max\{|a-\alpha|,
|b-\beta|\} \le \sigma < 0.1$. Then for $i,j \geq 0$,
$$|a^ib^j-\alpha^i\beta^j|  \le (1+2\sigma)\sigma.$$
\end{lemma}

\begin{proof}
We begin by bounding $|\alpha^i-a^i|$.  For $i=0$ or $\alpha \in
\{\pm1\}$, $\alpha^i=a^i$, so we can assume that $i>0$ and $\alpha \in
\{0,\pm 1/4, \pm 1/2\}$. Then
$$
\alpha^i - a^i = (\alpha-a) \sum_{k=0}^{i-1} \alpha^k a^{i-1-k}.
$$
If we apply the triangle inequality (together with $|\alpha-a| \le
\sigma$ and $|\alpha|, |a| \le 1/2+\sigma$), we find that
\begin{equation} \label{eq:alphaabd}
|\alpha^i - a^i| \le i(1/2+\sigma)^{i-1} \sigma.
\end{equation}
The function $n(1/2+\sigma)^{n-1}$ on nonnegative integers attains its
maximum value when $n = 2$, from which it follows that $|\alpha^i-a^i|
\le (1+2\sigma) \sigma$.  Note that this is the special case of the
lemma in which $\beta = b = \pm 1$.  Thus, we can henceforth assume
that neither $\beta$ nor $b$ is $\pm 1$, and by symmetry we can assume
the same for $\alpha$ and $a$.

For the general case, we notice that
$$
\alpha^i \beta^j - a^i b^j = (\alpha^i - a^i) \beta^j +
a^i(\beta^j - b^j),
$$
so
$$
|\alpha^i \beta^j - a^i b^j| \le |\alpha^i - a^i| \cdot|\beta^j| +
|a^i|\cdot|\beta^j - b^j|.
$$
It now follows from \eqref{eq:alphaabd} that
$$
|\alpha^i \beta^j - a^i b^j| \le (i+j)(1/2+\sigma)^{i+j-1} \sigma.
$$
The right-hand side is at most $(1+2\sigma)\sigma$, as before.
\end{proof}

Now we need to make precise the errors in all three
approximations in
\begin{align*}
\sum_{\alpha,\beta \in S} \alpha^i \beta^j
\widetilde{P}_{\gamma}(\alpha,\beta) &\approx  \sum_{w \in
\mathcal{C}_\Lambda} g_{i,j}(w)\\
&\approx \frac{196560}{\vol(S^{23})} \int_{S^{23}}
g_{i,j}(z) \, dz\\
&\approx \frac{196560}{\vol(S^{23})} G_{i,j}(\gamma).
\end{align*}

The first follows from Lemma~\ref{lemma:approx}:
$$
\left|\sum_{\alpha,\beta \in S} \alpha^i \beta^j
\widetilde{P}_\gamma (\alpha,\beta) - \sum_{w \in
\mathcal{C}_\Lambda} g_{i,j}(w)\right| \leq \sum_{w \in
\mathcal{C}_\Lambda} (1+2\sigma) \sigma = 196560(1+2\sigma)\sigma,
$$
because $|\alpha^i\beta^j - a^i b^j| \leq (1+2\sigma)\sigma$ for
$\langle w, x\rangle = a \approx \alpha$ and $\langle w, y\rangle
= b \approx \beta$.

The second we have already estimated sufficiently well in
Lemma~\ref{nearlydesign}, in terms of $|g_{i,j}|_2$. Because
$|g_{i,j}(z)| \le 1$ for all $z$, we have
$$
|g_{i,j}|_2 \le \sqrt{\vol(S^{23})} = \sqrt{24\cdot
\frac{\pi^{12}}{12!}} = \frac{\pi^6}{\sqrt{19958400}}.
$$
It follows that
$$
\left| \sum_{w \in \mathcal{C}_\Lambda} g_{i,j}(w) -
\frac{196560}{\vol(S^{23})} \int_{S^{23}} g_{i,j}(z) \, dz
\right| \le \frac{\pi^6}{\sqrt{19958400}} \cdot 2.50193 \cdot
10^{-5} < 5.3841 \cdot 10^{-6} .
$$

Finally, one can check by straightforward computation of each
case that if $\langle x,y \rangle$ differs from one of
$-1/2,-1/4,0,1/4,1/2$ by at most $\sigma$ (where $0 \le \sigma
\leq 1)$, and $0 \le i,j \le 4$, then
$$
\frac{196560}{\vol(S^{23})} \int_{S^{23}} g_{i,j}(z) \, dz
$$
differs by at most $8190\sigma$ from what it would be if $\sigma$
were zero. (In fact, if one expands this quantity as a power
series in $\sigma$, then the sum of the absolute values of the
coefficients is at most $8190$.)  Therefore the error in the last
approximation is at most $8190\sigma$.

Thus
$$
\sum_{\alpha,\beta \in S} \alpha^i \beta^j \widetilde{P}_\gamma
(\alpha,\beta) = \frac{196560}{\vol(S^{23})} G_{i,j}(\gamma) +
D_{i,j},
$$
where
$$
|D_{i,j}| < 196560(1+2\sigma)\sigma+ 8190\sigma + 5.3841 \cdot
10^{-6} < 6.7023\cdot 10^{-6}.
$$
As before the values of $\tilde{P}_\gamma(\pm 1,\alpha)$ and
$\tilde{P}_\gamma(\alpha, \pm 1)$ are known and are the same as
the corresponding $P_\gamma(\pm 1,\alpha)$ and
$P_\gamma(\alpha,\pm 1)$. Thus they serve as constants in the
equation and do not contribute to the error.

Let $A$ be the matrix of coefficients for these equations.  One
can check that $|A^{-1}|_\infty = 7225$.  (Here, $|\cdot|_\infty$
denotes the $\infty$-norm on matrices, which is induced by the
$\ell_\infty$ norm on vectors.  It is the maximum over all rows
of the sum of the absolute values of the elements in that row.) It
follows that the intersection numbers in $\mathcal{C}_\Lambda$
differ by at most
$$
7225 \cdot 6.7023\cdot 10^{-6} < 0.05
$$
from those in $\mathcal{C}_{24}$. Because they must be integers,
this proves that they are the same as in $\mathcal{C}_{24}$ (in
particular, they do not depend on the choice of $x$ and $y$,
although what we have proved is far stronger).

The computer file \texttt{verifyrest.txt} carries out all these
calculations (as well as those from several other points in this
paper).  In it, we assume without loss of generality that
$\gamma>0$ because $\widetilde{P}_\gamma(\alpha,\beta) =
\widetilde{P}_{-\gamma}(\alpha,-\beta)$.

\subsection{Uniqueness of association scheme}

The spherical code $\mathcal{C}_{24}$ determines a $6$-class association
scheme $\mathcal{A}_{24}$ if we partition elements $(x,y)$ of
$\mathcal{C}_{24} \times \mathcal{C}_{24}$ with $x\ne y$ according to their
inner products $\langle x,y\rangle$. We can similarly form the association
scheme $\mathcal{A}_\Lambda$ of the spherical code $\mathcal{C}_\Lambda$
coming from $\Lambda$, where this time we group elements according to their
approximate inner products. We wish to show that these two association
schemes are isomorphic. We will need to know that this $6$-class association
scheme $\mathcal{A}_{24}$ is uniquely determined by its size, valencies, and
intersection numbers. That can be proved as follows.

Let $N=196560$, let $C = 6/N = 1/32760$, and let $u_1,\dots,u_N$ be the
minimal vectors of the Leech lattice.  The following lemma restates the
(known) fact that $\Lambda_{24}$ is strongly eutactic.  We provide a
proof here to make this article more self-contained.

\begin{lemma}
\label{lemma:eutactic} For every $x \in \R^{24}$,
$$
\langle x,x\rangle = C\sum_{i=1}^N \langle x,u_i\rangle^2.
$$
\end{lemma}

\begin{proof}
Because the minimal vectors of $\Lambda_{24}$ form a spherical
$2$-design, the polynomial $y \mapsto \langle x,y \rangle^2$ has the
same average over $\{u_1,\dots,u_N\}$ and the sphere of radius $2$ in
$\R^{24}$, which is $4$ times the average over the unit sphere.  To
average over the unit sphere, it will be convenient to work with
orthonormal bases of $\R^{24}$.  For each orthonormal basis
$e_1,\dots,e_{24}$, we have $|x|^2 = \sum_i \langle x, e_i\rangle^2$.
If we average over all orthonormal bases, then each of
$e_1,\dots,e_{24}$ is uniformly distributed over the unit sphere, and
therefore the average of $y \mapsto \langle x,y \rangle^2$ over the
unit sphere is $|x|^2/24$. It follows that the average over the sphere
of radius $2$ is $|x|^2/6$, so
$$
\frac{1}{N}\sum_{i=1}^N \langle x, u_i \rangle^2 = |x|^2/6,
$$
as desired.
\end{proof}

\begin{theorem}\label{assoc_uniqueness}
There is only one $6$-class association scheme with the same size, valencies
and intersection numbers as the association scheme of minimal vectors of the
Leech lattice.
\end{theorem}

\begin{proof}
Let $\mathcal{A} = \{a_i : 1 \le i \le 196560\}$ be such an
association scheme, with $\mathcal{A}^2$ partitioned into classes
as
$$
\mathcal{A}^2 = \mathcal{A}_{1} \cup \mathcal{A}_{1/2} \cup
\mathcal{A}_{1/4} \cup \mathcal{A}_{0} \cup \mathcal{A}_{-1/4}
\cup \mathcal{A}_{-1/2} \cup \mathcal{A}_{-1}
$$
(labeled according to the corresponding inner products in the
Leech case, when the minimal vectors are rescaled to lie on the
unit sphere).

Let $A_1, A_{1/2}, \dots, A_{-1}$ denote the adjacency matrices of
the identity relation and the six classes; in other words,
$$
(A_\alpha)_{i,j} =
\begin{cases}
1 & \textup{if $(a_i,a_j) \in \mathcal{A}_\alpha,$ and}\\
0 & \textup{otherwise}.
\end{cases}
$$ These matrices are symmetric and commute with each other.
Their span forms an algebra, called the Bose-Mesner algebra,
whose product is determined by the valencies and intersection
numbers. (See \cite[p.~755]{BH}.)  Namely,
$$
A_{\alpha} A_{\beta} = \sum_{\gamma} P_{\gamma}(\alpha,\beta)
A_{\gamma}.
$$
Furthermore, note that because the $A_\alpha$'s have only $0$ or
$1$ as entries, are not identically zero, and sum to the all $1$'s
matrix, they must be linearly independent.  Thus, the Bose-Mesner
algebra is completely determined by the valencies and
intersection numbers, with no additional relations possible.

Let
$$
P = 4C \left(A_1 + \frac{1}{2} A_{1/2} +\frac{1}{4}A_{1/4}-
\frac{1}{4} A_{-1/4} -\frac{1}{2}A_{-1/2} - A_{-1}\right).
$$
If the association scheme $\mathcal{A}$ comes from the Leech
lattice, then $P$ is $C$ times the Gram matrix of the $196560$
minimal vectors (not rescaled to the unit sphere). We claim that
$P$ is a projection matrix, i.e., $P^2 = P$. One way to check
that is to use the valencies and intersection numbers to compute
$P^2$ and verify that it equals $P$.  That is somewhat cumbersome
to check by hand, so we will instead give a longer but
conceptually simpler proof.

First, we check it in the case of the actual Leech lattice
association scheme, using Lemma~\ref{lemma:eutactic}.  From
$$
\langle x,x\rangle = C \sum_{i=1}^N \langle x,u_i\rangle^2
$$
it follows that
\begin{align*}
\langle x,y \rangle &= \frac{1}{2}(\langle x+y,x+y \rangle -
\langle x,x \rangle - \langle y,y \rangle)\\
&= \frac{C}{2} \sum_{i=1}^N \left(\langle x+y,u_i \rangle^2 -
\langle x,u_i \rangle^2 - \langle y,u_i \rangle^2\right)\\
&= C \sum_{i=1}^N \langle x,u_i \rangle \langle y,u_i \rangle.
\end{align*}
Therefore the $(i,j)$ entry of $P^2$ is
$$
(P^2)_{i,j} = C^2 \sum_{k=1}^N \langle u_i,u_k\rangle \langle
u_j,u_k\rangle = C \langle u_i,u_j\rangle = P_{i,j},
$$
so for the Leech lattice, $P$ is a projection matrix.  Because the structure
of the Bose-Mesner algebra is completely determined by the valencies and
intersection numbers, the same must always be true.

Thus, $P$ is always a projection matrix (in fact, an orthogonal
projection because $P$ is symmetric).  The trace of $P$ is $4C$
times the trace of $A_1$, because no other $A_\alpha$'s have
entries on the diagonal.  The trace is therefore $4NC = 24$, so
$P$ projects onto a $24$-dimensional subspace.

Consider the images of the $N = 196560$ unit vectors (namely
$e_1,\dots,e_{196560}$) under $P$.  Their inner products are simply the
entries of $P$, since $\langle Pe_i,Pe_j\rangle = \langle e_i, P^2
e_j\rangle = \langle e_i, Pe_j\rangle = P_{i,j}$, and if one rescales
the vectors by $1/\sqrt{4C}$ so that they lie on the unit sphere, then
the result is a $196560$-point kissing configuration in $\R^{24}$.
(Appendix~\ref{appendix:background} reviews the kissing problem.) The
only such configuration is the kissing configuration of the Leech
lattice (see \cite{BS}, which was reprinted as Chapter~14 of
\cite{CS}), up to orthogonal transformations of $\R^{24}$. Thus, $P/C$
must be the Gram matrix of the minimal vectors in the Leech lattice.

It follows that $\mathcal{A}$ is isomorphic to $\mathcal{A}_{24}$
(in particular, $a_i \mapsto P(e_i)/\sqrt{4C}$ yields an
isomorphism). Thus, $\mathcal{A}_{24}$ is determined by its size
and intersection numbers, as desired.
\end{proof}

\section{Inner product bounds}\label{section:innerbounds}

We will first use the intersection numbers and isomorphism of association
schemes to prove better bounds on $\sigma$. We already know that for all
nearly minimal vectors $u$,
$$
4 \leq \langle u,u \rangle \leq 4(1+\varepsilon)^2 <
4+9\varepsilon.
$$

\begin{lemma}
For nearly minimal vectors $u,v$ with $\langle u,v \rangle \approx 2$,
$$
2-5\varepsilon \leq \langle u,v \rangle \leq 2+9 \varepsilon.
$$
\end{lemma}

\begin{proof}
We know that $u-v$ is nearly minimal, so
$$
\langle u-v,u-v \rangle \leq 4(1+\varepsilon)^2.
$$
It follows that
\begin{align*}
2\langle u,v \rangle &\geq - 4(1+2\varepsilon+\varepsilon^2) +
\langle u,u \rangle + \langle v,v \rangle \\
&\geq 8 - 4(1+2\varepsilon +\varepsilon^2)\\
&\geq 4-10 \varepsilon,
\end{align*}
which gives us one of the inequalities. Similarly,
$$
\langle u-v,u-v \rangle \geq 4
$$
gives us
\begin{align*}
2\langle u,v \rangle &\leq - 4 + \langle u,u \rangle + \langle
v,v \rangle \\ &\leq -4 + 8(1+2\varepsilon +\varepsilon^2)\\
&\leq 4+ 18\varepsilon,
\end{align*}
which is the other inequality.
\end{proof}

Note that these inequalities could be made slightly sharper, but
we prefer simpler numbers.

\begin{lemma}
For nearly minimal vectors $u,v$ with $\langle u,v \rangle \approx 0$,
$$
-14\varepsilon \leq \langle u,v \rangle \leq 14 \varepsilon.
$$
\end{lemma}

\begin{proof}
Since $P_0(1/2,1/2) = 44 \neq 0$, we can find $w$ with $\langle
u,w \rangle \approx 2$ and $\langle v,w \rangle \approx 2$. Then
$v-w$ is nearly minimal and we know from the previous lemma that
$$
2 - 5\varepsilon \leq \langle u,w \rangle \leq 2+9 \varepsilon
$$
and
$$
-2 -9 \varepsilon \leq \langle u,v-w \rangle \leq -2+5\varepsilon.
$$
Adding these inequalities gives us the result.
\end{proof}

\begin{lemma}
\label{lemma:leechconfig} If $u,v \in \Lambda_{24}$ are minimal
vectors with $\langle u,v \rangle = 1$, then there are minimal
vectors $w_1$, $w_2$, and $w_3$ satisfying $\langle u,w_i \rangle
= 2$, $\langle v,w_i \rangle = 0$, and $\langle w_i,w_j \rangle =
0$ for $i \ne j$.
\end{lemma}

\begin{proof}
We will use the fact that the automorphism group $\Co_0$ of the
Leech lattice acts transitively on pairs of minimal vectors with
a fixed inner product between them (see Theorem~3.13 in \cite{T}
for a proof). Thus it suffices to consider the case of a
particular pair $(u,v)$ of minimal vectors with inner product $1$.
Let
\begin{align*}
u &= \frac{1}{\sqrt 8}(1,1,\ldots,1,1,-3), \\
v &= \frac{1}{\sqrt 8}(0,0,\ldots,0,-4,-4), \\
w_1 &= \frac{1}{\sqrt 8}(2,2,2,2,2,2,2,2,0,0,\ldots,0), \\
w_2 &= \frac{1}{\sqrt 8}(0,0,\ldots,0,4,-4), \textup{ and} \\
w_3 &= \frac{1}{\sqrt
8}(0,0,0,0,0,0,0,0,2,2,0,0,2,2,0,0,2,2,0,0,2,2,0,0).
\end{align*}
(see \cite[p.~131]{CS} for a description of the minimal vectors).
It is easily checked that the inner products are as desired.
\end{proof}

\begin{lemma}\label{lemma:quarter}
For nearly minimal vectors $u,v$ with $\langle u,v \rangle \approx 1$,
$$
1-72\varepsilon \leq \langle u,v \rangle \leq 1+75\varepsilon.
$$
\end{lemma}

\begin{proof}
We know that the association schemes of the Leech lattice
$\Lambda_{24}$ and our given lattice $\Lambda$ are the same. Let
$u',v'$ be the corresponding vectors in the Leech lattice
(corresponding via some fixed isomorphism of association schemes).
Since $\langle u',v' \rangle = 1$, we know that we can find
$w_1'$, $w_2'$, and $w_3'$ in the Leech lattice satisfying
$\langle u',w_i' \rangle = 2$, $\langle v',w_i' \rangle = 0$, and
$\langle w_i', w_j' \rangle = 0$ by Lemma~\ref{lemma:leechconfig}.
Let $w_1,w_2,w_3$ be the corresponding vectors in $\Lambda$. Then
the relations $\langle u,w_i \rangle \approx 2$, $\langle v,w_i
\rangle \approx 0$, and $\langle w_i, w_j \rangle \approx 0$ must
hold in $\Lambda$. It follows by a short computation that
$2u-v-w_1-w_2-w_3$ is a nearly minimal vector. Therefore
\begin{align*}
4 \leq {}& \langle 2u-v-w_1-w_2-w_3, 2u-v-w_1-w_2-w_3\rangle\\
= {}& 4\langle u,u \rangle + \langle v,v\rangle + \sum_i \langle
w_i,w_i \rangle - 4 \langle u,v \rangle \\
&- 4\sum_i \langle u,w_i \rangle + 2\sum_i \langle
v,w_i \rangle + 2\sum_{i < j} \langle w_i, w_j \rangle.
\end{align*}
It follows that
\begin{align*}
4 \langle u,v \rangle \leq {}& {-4}+ 4\langle u,u \rangle + \langle
v,v\rangle + \sum_i \langle w_i,w_i \rangle\\
& - 4\sum_i \langle u,w_i \rangle + 2\sum_i \langle
v,w_i \rangle + 2\sum_{i < j} \langle w_i, w_j \rangle\\
\leq {}& {-4} + 8(4+9\varepsilon) -12(2-5\varepsilon) + 12(14
\varepsilon).
\end{align*}
Thus, $\langle u,v \rangle \leq 1+75 \varepsilon.$ Similarly,
\begin{align*}
4+9\varepsilon \geq {}& \langle 2u-v-w_1-w_2-w_3,
2u-v-w_1-w_2-w_3\rangle\\
= {}& 4\langle u,u \rangle + \langle v,v\rangle + \sum_i \langle
w_i,w_i \rangle - 4 \langle u,v \rangle\\
& - 4\sum_i \langle u,w_i \rangle + 2\sum_i \langle
v,w_i \rangle + 2\sum_{i < j} \langle w_i, w_j \rangle.
\end{align*}
It follows that
\begin{align*}
4 \langle u,v \rangle \geq {}& {-4}-9\varepsilon + 4\langle u,u
\rangle + \langle v,v\rangle + \sum_i \langle w_i,w_i \rangle \\
& - 4\sum_i \langle u,w_i \rangle
+ 2\sum_i \langle v,w_i \rangle + 2\sum_{i < j} \langle w_i, w_j \rangle \\
\geq {} & {-4}-9\varepsilon + 8(4) -12(2+9\varepsilon) + 12(-14
\varepsilon).
\end{align*}
Thus, $\langle u,v \rangle \geq 1-(285/4)\varepsilon \geq
1-72\varepsilon$.
\end{proof}

We have proved the following proposition.

\begin{proposition}\label{dotbounds}
If $u,v \in \Lambda$ are nearly minimal vectors, then $\langle
u,v\rangle$ differs from an element of $\{0, \pm 1, \pm 2, \pm
4\}$ by at most $75 \varepsilon$.
\end{proposition}

\section{A basis of nearly minimal vectors}\label{section:basis}

We wish to prove that $\Lambda$ must have a basis of nearly minimal
vectors.  We first prove that the nearly minimal vectors span
$\Lambda$. Let $x \in \Lambda$ be as small a vector as possible without
being in the span of the nearly minimal vectors.  Then there does not
exist a nearly minimal vector $u$ such that $|u-x| < |x|$. We know that
$|x| > {\sqrt {6} (1-\mu)}$ and $2 \le |u| \le 2(1+\varepsilon)$ for
each nearly minimal $u$. Because $|u-x| \ge |x|$, we have
$$
\langle u,x\rangle \le |u|^2/2 \le 2(1+\varepsilon)^2.
$$
Consider the unit vectors $x/|x|$ and $u/|u|$. We have
$$
\left\langle \frac{u}{|u|},
\frac{x}{|x|} \right\rangle \le
\frac{2(1+\varepsilon)^2}{{\sqrt{6}(1-\mu)} \cdot 2} <
\frac{1}{2}.
$$
If we extend the spherical code $\mathcal{C}_\Lambda$ of vectors of the
form $u/|u|$ to $\mathcal{C}_\Lambda \cup \{x/|x|\}$, then it will
contain $196561$ vectors without changing the minimal angle, and we
have seen that that is impossible.  Thus, the nearly minimal vectors do
span the lattice.

The same argument (with $\varepsilon=\mu=0$) also proves the
following fact: if the kissing configuration of a lattice is an
optimal kissing configuration for its dimension, then the lattice
is spanned by its minimal vectors.

Now let $B$ be $1/\sqrt{8}$ times the matrix
$$
\left(\textrm{
\begin{tiny}
\begin{tabular}{rrrrrrrrrrrrrrrrrrrrrrrr}
4 & {\hskip -9pt} $-4$ & 0 & 0 & 0 & 0 & 0 & 0 & 0 & 0 & 0 & 0 & 0 & 0 & 0 & 0 & 0 & 0 & 0 & 0 & 0 & 0 & 0 & 0\\
4 & 4 & 0 & 0 & 0 & 0 & 0 & 0 & 0 & 0 & 0 & 0 & 0 & 0 & 0 & 0 & 0 & 0 & 0 & 0 & 0 & 0 & 0 & 0\\
4 & 0 & 4 & 0 & 0 & 0 & 0 & 0 & 0 & 0 & 0 & 0 & 0 & 0 & 0 & 0 & 0 & 0 & 0 & 0 & 0 & 0 & 0 & 0\\
4 & 0 & 0 & 4 & 0 & 0 & 0 & 0 & 0 & 0 & 0 & 0 & 0 & 0 & 0 & 0 & 0 & 0 & 0 & 0 & 0 & 0 & 0 & 0\\
4 & 0 & 0 & 0 & 4 & 0 & 0 & 0 & 0 & 0 & 0 & 0 & 0 & 0 & 0 & 0 & 0 & 0 & 0 & 0 & 0 & 0 & 0 & 0\\
4 & 0 & 0 & 0 & 0 & 4 & 0 & 0 & 0 & 0 & 0 & 0 & 0 & 0 & 0 & 0 & 0 & 0 & 0 & 0 & 0 & 0 & 0 & 0\\
4 & 0 & 0 & 0 & 0 & 0 & 4 & 0 & 0 & 0 & 0 & 0 & 0 & 0 & 0 & 0 & 0 & 0 & 0 & 0 & 0 & 0 & 0 & 0\\
2 & 2 & 2 & 2 & 2 & 2 & 2 & 2 & 0 & 0 & 0 & 0 & 0 & 0 & 0 & 0 & 0 & 0 & 0 & 0 & 0 & 0 & 0 & 0\\
4 & 0 & 0 & 0 & 0 & 0 & 0 & 0 & 4 & 0 & 0 & 0 & 0 & 0 & 0 & 0 & 0 & 0 & 0 & 0 & 0 & 0 & 0 & 0\\
4 & 0 & 0 & 0 & 0 & 0 & 0 & 0 & 0 & 4 & 0 & 0 & 0 & 0 & 0 & 0 & 0 & 0 & 0 & 0 & 0 & 0 & 0 & 0\\
4 & 0 & 0 & 0 & 0 & 0 & 0 & 0 & 0 & 0 & 4 & 0 & 0 & 0 & 0 & 0 & 0 & 0 & 0 & 0 & 0 & 0 & 0 & 0\\
2 & 2 & 2 & 2 & 0 & 0 & 0 & 0 & 2 & 2 & 2 & 2 & 0 & 0 & 0 & 0 & 0 & 0 & 0 & 0 & 0 & 0 & 0 & 0\\
4 & 0 & 0 & 0 & 0 & 0 & 0 & 0 & 0 & 0 & 0 & 0 & 4 & 0 & 0 & 0 & 0 & 0 & 0 & 0 & 0 & 0 & 0 & 0\\
2 & 2 & 0 & 0 & 2 & 2 & 0 & 0 & 2 & 2 & 0 & 0 & 2 & 2 & 0 & 0 & 0 & 0 & 0 & 0 & 0 & 0 & 0 & 0\\
2 & 0 & 2 & 0 & 2 & 0 & 2 & 0 & 2 & 0 & 2 & 0 & 2 & 0 & 2 & 0 & 0 & 0 & 0 & 0 & 0 & 0 & 0 & 0\\
2 & 0 & 0 & 2 & 2 & 0 & 0 & 2 & 2 & 0 & 0 & 2 & 2 & 0 & 0 & 2 & 0 & 0 & 0 & 0 & 0 & 0 & 0 & 0\\
4 & 0 & 0 & 0 & 0 & 0 & 0 & 0 & 0 & 0 & 0 & 0 & 0 & 0 & 0 & 0 & 4 & 0 & 0 & 0 & 0 & 0 & 0 & 0\\
2 & 0 & 2 & 0 & 2 & 0 & 0 & 2 & 2 & 2 & 0 & 0 & 0 & 0 & 0 & 0 & 2 & 2 & 0 & 0 & 0 & 0 & 0 & 0\\
2 & 0 & 0 & 2 & 2 & 2 & 0 & 0 & 2 & 0 & 2 & 0 & 0 & 0 & 0 & 0 & 2 & 0 & 2 & 0 & 0 & 0 & 0 & 0\\
2 & 2 & 0 & 0 & 2 & 0 & 2 & 0 & 2 & 0 & 0 & 2 & 0 & 0 & 0 & 0 & 2 & 0 & 0 & 2 & 0 & 0 & 0 & 0\\
0 & 2 & 2 & 2 & 2 & 0 & 0 & 0 & 2 & 0 & 0 & 0 & 2 & 0 & 0 & 0 & 2 & 0 & 0 & 0 & 2 & 0 & 0 & 0\\
0 & 0 & 0 & 0 & 0 & 0 & 0 & 0 & 2 & 2 & 0 & 0 & 2 & 2 & 0 & 0 & 2 & 2 & 0 & 0 & 2 & 2 & 0 & 0\\
0 & 0 & 0 & 0 & 0 & 0 & 0 & 0 & 2 & 0 & 2 & 0 & 2 & 0 & 2 & 0 & 2 & 0 & 2 & 0 & 2 & 0 & 2 & 0\\
{\hskip -9pt}$-3$ & 1 & 1 & 1 & 1 & 1 & 1 & 1 & 1 & 1 & 1 & 1 & 1
& 1 & 1 & 1 & 1 & 1 & 1 & 1 & 1 & 1 & 1 & 1
\end{tabular}
\end{tiny}}
\right).
$$

The rows of $B$ form a basis of $\Lambda_{24}$ consisting of minimal
vectors (see Figure~4.12 in \cite{CS}).  One can compute the inverse
matrix and check that its entries are integers divided by $\sqrt{8}$.
The largest entry of $B^{-1}$, in absolute value, is $-13/\sqrt{8}$, so
$|B^{-1}|_\infty \le 24 \cdot 13/\sqrt{8}$.

Let $u_1',\ldots,u_{24}'$ be the rows of $B$. Every minimal
vector $u'$ is a linear combination $\sum_i c_i u_i'$ of this
basis, and the coefficients are bounded by
$$
|c_i| \leq |B^{-1}|_\infty |u'|_\infty \leq (24\cdot 13/\sqrt{8})
\cdot ({4}/{\sqrt 8}) = 156.
$$
Consider the corresponding nearly minimal vectors
$u_1,\dots,u_{24}$ of $\Lambda$ (under an isomorphism of
association schemes).   We next prove that they form a basis.

We need only check that $u_1,\dots,u_{24}$ span the nearly minimal
vectors of $\Lambda$. If $u$ is any nearly minimal vector, the
isomorphism of association schemes gives us numbers
$c_1,\ldots,c_{24}$ such that $u$ should be $\sum_i c_i u_i$
(i.e., the corresponding equality is true for the Leech lattice).
We first check that $\sum_i c_i u_i$ is a nearly minimal vector,
and then that it equals $u$. To check that it is nearly minimal,
we just need to know that all inner products of nearly minimal
vectors are within $75\varepsilon$ of what they are in the Leech
lattice (Proposition~\ref{dotbounds}). When we compute the norm of
$\sum_i c_i u_i$, each inner product $\langle u_i, u_j \rangle$
could be off by as much as $75 \varepsilon$, and multiplied by up
to $156^2$. There are $24^2$ such pairs, for a total error of at
most $156^2 \cdot 24^2 \cdot 75\varepsilon = 1051315200
\varepsilon$, which is minuscule (less than $10^{-17}$). Because
the next smallest vectors beyond the nearly minimal vectors have
norms at least $\sqrt{6}(1-\mu)$, we conclude that $\sum_i c_i
u_i$ must be nearly minimal.

To check that $u = \sum_i c_i u_i$, we need only compute their
inner product and verify that it is approximately $4$. This time,
the maximum error is $24 \cdot 156 \cdot 75\varepsilon$, which is
again small enough by a huge margin.

Thus, $u_1,\dots,u_{24}$ form a basis of $\Lambda$, and their
inner products are all within $75\varepsilon$ of what they would
be in the Leech lattice.

\section{Local optimality of $\Lambda_{24}$}
\label{section:localopt}

In this section, we will prove in detail that the Leech lattice is
locally optimal, and provide quantitative bounds.  We will follow
the notation and techniques from \cite{GL} closely.  The basic
result is Voronoi's theorem, which says that a lattice is a local
maximum for sphere packing density if and only if it is perfect
and eutactic.  In this section of the paper only, $\Lambda$ will
denote an arbitrary lattice in $\R^n$.

It is convenient to work in terms of the quadratic form $Q$ associated
to $\Lambda$.  Choose a lattice basis $\{ b_i \}$, and for a vector $x
\in \R^n$ (with coordinates $x_1,\dots,x_n$) define
$$
Q(x) = \left\langle \sum_{i=1}^n x_i b_i, \sum_{i=1}^n x_i b_i
\right\rangle.
$$
The matrix of $Q$ is $S = (s_{i,j})_{i,j=1}^n$, where $s_{i,j} =
\langle b_i, b_j \rangle$.

Let $M$ denote the minimal nonzero norm of $\Lambda$, and let
$u_1, \dots, u_N \in \Z^n$  be the coefficient vectors of the
minimal vectors in terms of the basis $\{ b_i \}$.  Thus, for $1
\le i \le N$,
$$
Q(u_i) = M.
$$
Recall that $Q$ is \textit{perfect\/} if these equations
completely determine the quadratic form $Q$.  Equivalently, every
quadratic form that vanishes at $u_1,\dots,u_N$ must vanish
everywhere.

Let $D$ denote the determinant of $S$, and let $\tilde{S} =
(\tilde{s}_{i,j})_{i,j=1}^n$ denote the adjoint matrix, where
$$
\tilde{s}_{i,j} = \frac{\partial D}{\partial s_{i,j}}.
$$
Strictly speaking this is an abuse of notation, but of course it
means we take the partial derivatives of $D$ as if the entries of
$S$ were variables, and then substitute their actual values.

In other words, $S \tilde{S}$ is $D$ times the identity matrix.
(It might seem that a transpose is missing, but note that all our
matrices are symmetric.)  Let $\tilde{Q}$ be the quadratic form
with matrix $\tilde{S}$.  Then $\Lambda$ is \textit{eutactic\/} if
there are positive numbers $d_1,\dots,d_N$ such that for all $x
\in \R^n$,
$$
\tilde{Q}(x) = \sum_{k=1}^N d_k \langle u_k, x \rangle^2.
$$

It is known that $\Lambda_{24}$ is perfect and eutactic, with
$d_1 = \dots = d_{196560} = 1/32760$. For completeness we sketch a
proof.  Lemma~\ref{lemma:eutactic} proves that the Leech lattice
is eutactic.

\begin{lemma}
The Leech lattice $\Lambda_{24}$ is perfect.
\end{lemma}

\begin{proof}
Suppose $Q$ is a quadratic form that vanishes on the minimal
vectors. Then the symmetric bilinear form $B$ corresponding to
$Q$ satisfies $B(u_i,u_i) = 0$. We first show that $B(u,v) = 0$
for all minimal vectors $u,v$. If $\langle u,v \rangle = 2$ we
use the fact that $u -v$ is minimal to see that
$$B(u,v) = \frac{1}{2}\left(B(u,u)+B(v,v)-B(u-v,u-v) \right) = 0.
$$
If $\langle u,v \rangle = 0$, then since $P_0(1/2,-1/2) = 44 \neq
0$, we can find $w$ with $\langle u,w \rangle = 2$ and $\langle
v,w \rangle = -2$. Then $v+w$ is nearly minimal and we know
$$
0 = B(u,w) = B(u,v+w);
$$
subtracting these two gives us the result. If $\langle u,v \rangle
= 1$, then by Lemma~\ref{lemma:leechconfig} there are minimal
vectors $w_1$, $w_2$, and $w_3$ with $\langle u,w_i \rangle = 2$,
$\langle v,w_i \rangle = 0$, and $\langle w_i, w_j \rangle = 0$
for $i \ne j$. It follows that $2u-v-w_1-w_2-w_3$ is a minimal
vector. Therefore
\begin{align*}
0 = {}&  B( 2u-v-w_1-w_2-w_3, 2u-v-w_1-w_2-w_3) \\ = {}& 4 B(u,u)+
B(v,v) + \sum_i B(w_i,w_i)-  4B(u,v)\\
&  - 4\sum_i B(u,w_i) + 2\sum_i B(v,w_i)+ 2\sum_{i < j}
B(w_i, w_j) \\
= {}& {- 4}B(u,v) + 0.
\end{align*}
This forces $B(u,v) = 0$. Finally, the result for $\langle u,v \rangle < 0$
follows from the above because $B(u,v) = -B(u,-v)$ and $\langle u, -v \rangle
> 0$.

To conclude the proof, we use the above information on a basis of minimal
vectors to see that $B$ is identically zero, and hence $Q$ is as well.
\end{proof}

We begin by proving one direction of Voronoi's theorem. This proof is the one
given in \cite[\S 39]{GL}, where one can also find a proof of the converse.
We give the details of this direction because we will need to examine it in
detail to derive quantitative estimates, and because it is the only direction
needed here.

\begin{theorem}[Voronoi \cite{Vo}]\label{theorem:voronoi}
If $\Lambda$ is perfect and eutactic, then it is a strict local
maximum for density.
\end{theorem}

\begin{proof}
We wish to show that if $\Lambda$ is perturbed slightly (other than
simply by scaling and isometries), then $D^{-1/n}M$ must strictly
decrease. We perturb $S=(s_{i,j})$ by changing it to $(s_{i,j}+\rho
t_{i,j})$, where $t_{j,i}=t_{i,j}$ and $\rho>0$ is small.  We assume
$(t_{i,j})$ is not identically zero.  Let $Q_\rho$ denote the
corresponding quadratic form, and let $D_\rho$ be the determinant of
the matrix $(s_{i,j}+\rho t_{i,j})$.  We will show that for any fixed
matrix $(t_{i,j})$ not proportional to $S$, if $\rho$ is sufficiently
small, then there exists a $k$ such that
$$
D_\rho^{-1/n} Q_\rho(u_k) < D^{-1/n} M.
$$
Because $Q_\rho(u_k)$ is no smaller than the minimal norm of $Q_\rho$
(i.e., the minimum value of $Q_\rho$ on $\Z^n \setminus \{0\}$), this
inequality is what we want.

First, note that without loss of generality we can assume that
\begin{equation}
\label{zerosum}
\sum_{i,j} \tilde{s}_{i,j} t_{i,j} = 0.
\end{equation}
The reason is that
$$
\sum_{i,j} \tilde{s}_{i,j} s_{i,j} = nD \ne 0.
$$
Every perturbation can be broken up into the sum of a
perturbation proportional to $S$ and a perturbation satisfying
\eqref{zerosum}.  If we deal with the latter, then we can ignore
the former (which rescales the lattice but does not change its
packing density).

Thus, we assume \eqref{zerosum} from now on. Then the determinant
$D_\rho$ is given by
\begin{align*}
D_\rho = \det(s_{i,j}+ \rho t_{i,j}) &= \det(s_{i,j}) + \rho
\sum_{i,j}t_{i,j}\tilde{s}_{i,j} + O(\rho^2) \\
&= D + O(\rho^2) \quad \textup{as $\rho \rightarrow 0$}.
\end{align*}
Because $\Lambda$ is eutactic,
$$
\tilde{Q}(x) = \sum_{k=1}^N d_k \langle u_k, x \rangle^2
$$
for all $x \in \R^n$.  The associated symmetric bilinear form is
$$
\sum_{k=1}^N d_k \langle u_k, x \rangle\langle u_k, y \rangle,
$$
from which it follows that
$$
\tilde{s}_{i,j} = \sum_{k=1}^{N} d_k (u_{k})_i(u_{k})_j,
$$
where $(u_k)_i$ denotes the $i$-th coefficient of $u_k$. Therefore
$$
\sum_{k=1}^N d_k \sum_{i,j} t_{i,j} (u_{k})_i (u_{k})_j =
\sum_{i,j}\tilde{s}_{i,j} t_{i,j} = 0.
$$
Because $\Lambda$ is perfect, the inner sum on the left-hand side in the
above equation cannot vanish for all $k$.  Therefore, there exists a $k$ for
which it is negative, say
$$
\sum_{i,j} t_{i,j} (u_{k})_i (u_{k})_j \le -\alpha
$$
with $\alpha>0$. Then
$$
Q_\rho(u_k) = \sum (s_{i,j}+\rho t_{i,j}) (u_k)_i (u_k)_j \le M
-\rho\alpha,
$$
and hence
$$
D_\rho^{-1/n} Q_\rho(u_k) \le D^{-1/n} M (1-\rho\alpha/M)
(1+O(\rho^2))^{-1/n} < D^{-1/n}M
$$
if $\rho$ is positive and small enough.  This proves that
$\Lambda$ is a strict local optimum for density when $(s_{i,j})$
is perturbed in the direction of $(t_{i,j})$.  In fact, the
choices of $\alpha$ and the implicit constant in the big-$O$ can
be made uniformly in $(t_{i,j})$, given $\sum_{i,j}
|t_{i,j}|^2=1$. Thus $\Lambda$ is a strict local optimum for
density.
\end{proof}

We now compute numerical bounds on the perturbations.  We use the
same basis of $\Lambda_{24}$ as before.  The corresponding Gram
matrix is
$$
\left(\textrm{\begin{tiny}
\begin{tabular}{rrrrrrrrrrrrrrrrrrrrrrrr}
4 & 0 & 2 & 2 & 2 & 2 & 2 & 0 & 2 & 2 & 2 & 0 & 2 & 0 & 1 & 1 & 2 & 1 & 1 & 0 & {\hskip -9pt} $-1$ & 0 & 0 & {\hskip -9pt} $-2$\\
0 & 4 & 2 & 2 & 2 & 2 & 2 & 2 & 2 & 2 & 2 & 2 & 2 & 2 & 1 & 1 & 2 & 1 & 1 & 2 & 1 & 0 & 0 & {\hskip -9pt} $-1$\\
2 & 2 & 4 & 2 & 2 & 2 & 2 & 2 & 2 & 2 & 2 & 2 & 2 & 1 & 2 & 1 & 2 & 2 & 1 & 1 & 1 & 0 & 0 & {\hskip -9pt} $-1$\\
2 & 2 & 2 & 4 & 2 & 2 & 2 & 2 & 2 & 2 & 2 & 2 & 2 & 1 & 1 & 2 & 2 & 1 & 2 & 1 & 1 & 0 & 0 & {\hskip -9pt} $-1$\\
2 & 2 & 2 & 2 & 4 & 2 & 2 & 2 & 2 & 2 & 2 & 1 & 2 & 2 & 2 & 2 & 2 & 2 & 2 & 2 & 1 & 0 & 0 & {\hskip -9pt} $-1$\\
2 & 2 & 2 & 2 & 2 & 4 & 2 & 2 & 2 & 2 & 2 & 1 & 2 & 2 & 1 & 1 & 2 & 1 & 2 & 1 & 0 & 0 & 0 & {\hskip -9pt} $-1$\\
2 & 2 & 2 & 2 & 2 & 2 & 4 & 2 & 2 & 2 & 2 & 1 & 2 & 1 & 2 & 1 & 2 & 1 & 1 & 2 & 0 & 0 & 0 & {\hskip -9pt} $-1$\\
0 & 2 & 2 & 2 & 2 & 2 & 2 & 4 & 1 & 1 & 1 & 2 & 1 & 2 & 2 & 2 & 1 & 2 & 2 & 2 & 2 & 0 & 0 & 1\\
2 & 2 & 2 & 2 & 2 & 2 & 2 & 1 & 4 & 2 & 2 & 2 & 2 & 2 & 2 & 2 & 2 & 2 & 2 & 2 & 1 & 1 & 1 & {\hskip -9pt} $-1$\\
2 & 2 & 2 & 2 & 2 & 2 & 2 & 1 & 2 & 4 & 2 & 2 & 2 & 2 & 1 & 1 & 2 & 2 & 1 & 1 & 0 & 1 & 0 & {\hskip -9pt} $-1$\\
2 & 2 & 2 & 2 & 2 & 2 & 2 & 1 & 2 & 2 & 4 & 2 & 2 & 1 & 2 & 1 & 2 & 1 & 2 & 1 & 0 & 0 & 1 & {\hskip -9pt} $-1$\\
0 & 2 & 2 & 2 & 1 & 1 & 1 & 2 & 2 & 2 & 2 & 4 & 1 & 2 & 2 & 2 & 1 & 2 & 2 & 2 & 2 & 1 & 1 & 1\\
2 & 2 & 2 & 2 & 2 & 2 & 2 & 1 & 2 & 2 & 2 & 1 & 4 & 2 & 2 & 2 & 2 & 1 & 1 & 1 & 1 & 1 & 1 & {\hskip -9pt} $-1$\\
0 & 2 & 1 & 1 & 2 & 2 & 1 & 2 & 2 & 2 & 1 & 2 & 2 & 4 & 2 & 2 & 1 & 2 & 2 & 2 & 2 & 2 & 1 & 1\\
1 & 1 & 2 & 1 & 2 & 1 & 2 & 2 & 2 & 1 & 2 & 2 & 2 & 2 & 4 & 2 & 1 & 2 & 2 & 2 & 2 & 1 & 2 & 1\\
1 & 1 & 1 & 2 & 2 & 1 & 1 & 2 & 2 & 1 & 1 & 2 & 2 & 2 & 2 & 4 & 1 & 2 & 2 & 2 & 2 & 1 & 1 & 1\\
2 & 2 & 2 & 2 & 2 & 2 & 2 & 1 & 2 & 2 & 2 & 1 & 2 & 1 & 1 & 1 & 4 & 2 & 2 & 2 & 1 & 1 & 1 & {\hskip -9pt} $-1$\\
1 & 1 & 2 & 1 & 2 & 1 & 1 & 2 & 2 & 2 & 1 & 2 & 1 & 2 & 2 & 2 & 2 & 4 & 2 & 2 & 2 & 2 & 1 & 1\\
1 & 1 & 1 & 2 & 2 & 2 & 1 & 2 & 2 & 1 & 2 & 2 & 1 & 2 & 2 & 2 & 2 & 2 & 4 & 2 & 2 & 1 & 2 & 1\\
0 & 2 & 1 & 1 & 2 & 1 & 2 & 2 & 2 & 1 & 1 & 2 & 1 & 2 & 2 & 2 & 2 & 2 & 2 & 4 & 2 & 1 & 1 & 1\\
{\hskip -9pt} $-1$ & 1 & 1 & 1 & 1 & 0 & 0 & 2 & 1 & 0 & 0 & 2 & 1 & 2 & 2 & 2 & 1 & 2 & 2 & 2 & 4 & 2 & 2 & 2\\
0 & 0 & 0 & 0 & 0 & 0 & 0 & 0 & 1 & 1 & 0 & 1 & 1 & 2 & 1 & 1 & 1 & 2 & 1 & 1 & 2 & 4 & 2 & 2\\
0 & 0 & 0 & 0 & 0 & 0 & 0 & 0 & 1 & 0 & 1 & 1 & 1 & 1 & 2 & 1 & 1 & 1 & 2 & 1 & 2 & 2 & 4 & 2\\
{\hskip -9pt} $-2$ & {\hskip -9pt} $-1$ & {\hskip -9pt} $-1$ &
{\hskip -9pt} $-1$ & {\hskip -9pt} $-1$ & {\hskip -9pt} $-1$ &
{\hskip -9pt} $-1$ & 1 & {\hskip -9pt} $-1$ & {\hskip -9pt} $-1$
& {\hskip -9pt} $-1$ & 1 & {\hskip -9pt} $-1$ & 1 & 1 & 1 &
{\hskip -9pt} $-1$ & 1 & 1 & 1 & 2 & 2 & 2 & 4
\end{tabular}
\end{tiny}}
\right).
$$

Suppose we assume that the Gram matrix entries are perturbed by $\rho
t_{i,j}$ where $\max\{|t_{i,j}|\} \leq 1$ (and \eqref{zerosum} holds---we
will eventually deal with the case where it does not). Then, with the same
notation as in the proof of Theorem~\ref{theorem:voronoi},
$$
D_\rho-D = \sum_{T} \rho^{\dim T} \det(T) \det(\tilde{T}),
$$
where $T$ ranges over all nonvacuous minors of $(t_{i,j})$ and
$\det(\tilde{T})$ is the corresponding cofactor of the matrix $S$.
(The $D$ term corresponds to the case where $T$ is vacuous, i.e.,
contains no rows or columns.)  This expansion follows from
combining the Laplace expansion (see \cite[\S 95]{Mu}) with
multilinearity, and is known as Albeggiani's theorem \cite[\S
96]{Mu}.

Now let $\dim T = k$. Then $|\det(T)| \leq k^{k/2}$ by Hadamard's inequality,
since the length of each row of $T$ is at most $\sqrt{k}$ (see (7.8.2) in
\cite{HJ}). It follows that the absolute value of the sum of $\rho^k \det(T)
\det(\tilde{T})$ over $k$-dimensional $T$ is bounded by $k^{k/2} \rho^k A_k$,
where $A_k$ is the sum of the absolute values of the $(24-k)$-dimensional
minors of the Gram matrix. For $k \geq 3$ we use the simple bound
$$
A_k \leq \binom{24}{k}^2 \cdot (4^2 +(24-k-1)(2^2))^{(24-k)/2}
$$
(the first factor is the number of $(24-k) \times (24-k)$ minors, and the
second is the above bound on the determinant of the cofactor, because the
largest entry in each row of $S$ is $4$ and the other entries are at most $2$
in absolute value). For $k=2$ we explicitly compute the sum of absolute
values of the $22 \times 22$ minors, and find that it is $818153$. Putting
all this together, we see that for $0 < \rho < 10^{-20}$,
\begin{equation}\label{eq:drho}
D_\rho \geq 1 - \rho^2\big(2^{2/2}\cdot 818153 + 2\cdot 10^8\big),
\end{equation}
where the $2\cdot 10^8$ term bounds the contribution from all
higher powers of $\rho$. (Recall that $D=1$ for the Leech
lattice.)  See the computer file \texttt{verifygram.txt} for the
details of this calculation.

Next, we find an $\alpha$ that works for any choice of $t_{i,j}$
such that $\max\{|t_{i,j}|\} = 1$. This is done by linear
programming as follows.

We find an $\alpha > 0$ such that for all $(i_0,j_0)$, and for all
$t_{i,j}$ subject to the constraints  $-1 \leq t_{i,j} \leq 1$ for
$(i,j)\neq (i_0,j_0)$, $t_{i,j}= t_{j,i}$, $t_{i_0,j_0} = 1$, and
$\sum_{i,j} \tilde{s}_{i,j} t_{i,j} = 0$, the following inequality
holds for some minimal vector $u_k$:
$$
\sum_{i,j}t_{i,j}(u_{k})_i(u_{k})_j \leq -\alpha.
$$
We do the same for $t_{i_0,j_0} = -1$. All these linear programs
can be solved by computer, and it appears that $\alpha=1/23$
works and is the largest possible value of $\alpha$. However,
that is the result of floating point calculations that are not
rigorous. We will be content with proving (without computer
assistance) that $\alpha = 4/1055$ satisfies these properties.
This weaker bound will suffice for our purposes and is proved in
Section~\ref{section:alpha}.

We conclude that for every perturbation by $t_{i,j}$ where
$\max\{|t_{i,j}|\} = \rho$ and $\sum_{i,j} \tilde{s}_{i,j}t_{i,j}
= 0$, there exists a $k$ such that for $0 < \rho < 10^{-20}$,
\begin{equation}\label{eq:upperbd}
\frac{D_\rho^{-1/24} Q_\rho(u_k)}{D^{-1/24} \cdot M} \le
(1-\rho/1055) (1-(2\cdot10^8 + 2\cdot818153)\rho^2)^{-1/24}.
\end{equation}
We used here that $M = 4$ for $\Lambda_{24}$. The upper bound in
\eqref{eq:upperbd} is strictly less than $1$ when $0 < \rho <
10^{-20}$. (Note that for notational convenience we have absorbed
the factor $\rho$ into the perturbations $t_{i,j}$.)

The last remaining issue is that our perturbation may not satisfy
$\sum_{i,j} \tilde{s}_{i,j}t_{i,j} = 0$.  Suppose our perturbed matrix
entries are $s_{i,j}+\eta_{i,j}$.  Let
$$
\Delta = \sum_{i,j} \tilde{s}_{i,j} \eta_{i,j}.
$$
By Proposition~\ref{dotbounds}, $|\eta_{i,j}| \le 75\varepsilon$.  It
follows that $|\Delta| \le 152100 \varepsilon$, since the sum of the
absolute values of the entries of $\tilde{S}$ is $2028$.

If we divide the quadratic form by $1+\Delta/24$, which is
nonzero, then it is equivalent to $\Lambda_{24}$ perturbed by
$$
t_{i,j} = \frac{s_{i,j} + \eta_{i,j}}{1+\Delta/24} - s_{i,j},
$$
where now
$$
\sum_{i,j} \tilde{s}_{i,j}t_{i,j} = 0,
$$
because
$$
\sum_{i,j} \tilde{s}_{i,j}s_{i,j} = 24D = 24.
$$
To conclude our proof, we need only check that $|t_{i,j}| <
10^{-20}$.  (Note that if $t_{i,j}=0$ for all $i,j$, then the
perturbed quadratic form is proportional to the original one and
therefore equal to it because they have the same determinant.)

Using $|\Delta| \le 152100 \varepsilon$, $|\eta_{i,j}| \le
75\varepsilon$, and $|s_{i,j}| \le 4$, we have
\begin{equation} \label{eq:final}
|t_{i,j}| =
\left|\frac{\eta_{i,j}-s_{i,j}\Delta/24}{1+\Delta/24}\right| \le
\frac{75\varepsilon + 4\cdot 152100
\varepsilon/24}{1-152100\varepsilon/24} < 1.8 \cdot 10^{-22}.
\end{equation}
Because $1.8 \cdot 10^{-22} < 10^{-20}$, we find that the dense
lattice $\Lambda$ is close enough to $\Lambda_{24}$ to conclude
that $\Lambda$ is either the same as $\Lambda_{24}$ (up to
isometries of $\R^{24}$) or strictly less dense than
$\Lambda_{24}$.

This completes the proof of our main theorem, except for the
postponed computation of $\alpha$ in Section~\ref{section:alpha}:

\begin{theorem}\label{theorem:main}
The Leech lattice is the unique densest lattice in $\R^{24}$, up
to scaling and isometries of $\R^{24}$.
\end{theorem}

Note that the reason why scaling ambiguity appears in the theorem
statement but not in the above proof is that we fixed
$|\Lambda|=1$.

\section{Computation of $\alpha$}
\label{section:alpha}

Suppose $i_0,j_0 \in \{1,2,\dots,24\}$ and $t = \pm 1$. We wish to find
a number $\alpha>0$ such that whenever  $-1 \leq t_{i,j} \leq 1$ for
$(i,j)\neq (i_0,j_0)$, $t_{i,j} = t_{j,i}$, $t_{i_0,j_0} = t$, and
$\sum_{i,j} \tilde{s}_{i,j} t_{i,j} = 0$, there is a $k$ such that
$$
\sum_{i,j}t_{i,j}(u_{k})_i(u_{k})_j \le -\alpha.
$$

We will see that we can take $\alpha=4/1055$, whatever $i_0$,
$j_0$, and $t$ are.

\begin{lemma}
\label{lemma:ortho} Each minimal vector of the Leech lattice is
contained in a set of twenty-four orthogonal minimal vectors.
\end{lemma}

In fact, more is true: the minimal vectors can be partitioned into
$4095$ sets $\{\pm v_1,\dots, \pm v_{24}\}$ with $v_i$ and $v_j$
orthogonal for $i \ne j$.  See footnote~3 of \cite[p.~6]{E1} for an
elegant proof.

\begin{proof}
Because the automorphism group of the Leech lattice acts
transitively on the minimal vectors, we need only verify that
there exists a set of twenty-four orthogonal minimal vectors. Let
$w_i$ be the vector
$$
\frac{1}{\sqrt 8}(0,\dots,0,4,4,0,\dots,0),
$$
where only coordinates $i$ and $i+1$ are nonzero, and let $v_i$
be the vector
$$
\frac{1}{\sqrt 8}(0,\dots,0,4,-4,0,\dots,0).
$$
These are all minimal vectors in the Leech lattice, and
$w_1,v_1,w_3,v_3,\dots,w_{23},v_{23}$ is an orthogonal basis of
$\R^{24}$.
\end{proof}

Let $T$ denote the matrix $(t_{i,j})$.  We will also write $T(v)
= \langle v, Tv \rangle$ and $T(u,v) = \langle u, Tv \rangle$. In
these terms, our goal is to show that given the assumptions on
$T$, there is some $k$ such that $T(u_k) \le -\alpha$.

We will make use of the following lemma, which depends only on the
hypotheses listed in its statement:

\begin{lemma}
\label{lemma:orthsum} Let $v_1,\dots,v_{24}$ be orthogonal minimal
vectors in the Leech lattice.  If $\sum_{i,j} \tilde{s}_{i,j}
t_{i,j} = 0$, then
$$
\sum_{i=1}^{24} T(v_i) = 0.
$$
\end{lemma}

\begin{proof}
Let $B$ be the matrix whose rows are the coordinates of
$v_1,\dots,v_{24}$ relative to the basis we have chosen for the
Leech lattice. Then we see by orthogonality that $BSB^t = 4I$.
Thus, $S = 4 B^{-1}(B^t)^{-1}$.  It follows that
$$
\Tr(BTB^t) = \Tr(B^tBT) = 4\Tr(S^{-1}T) = 4\sum_{i,j}
\tilde{s}_{i,j} t_{i,j} = 0,
$$
which implies
\[
\sum_{i=1}^{24} T(v_i) = 0. \qedhere
\]
\end{proof}

Let us rephrase our basic problem, and slightly weaken the
hypotheses, as follows. Suppose $\sum_{i,j} \tilde{s}_{i,j}
t_{i,j} = 0$, and we are given minimal vectors $x,y$ such that
$$
\langle x,y \rangle =  \beta \in \{\pm 4, \pm 2, \pm 1, \pm 0\}
$$
and $T(x,y) = t \ne 0$. We wish to find $\alpha>0$ such that there
must always exist a minimal vector $w$ with $T(w,w) \leq -\alpha
|t|$ under these hypotheses. We will apply it with $t = \pm 1$.
(Note that we are no longer assuming $|t_{i,j}| \le 1$.)

If we prove a bound for a certain $(\beta,t)$ we automatically get
the same bound of $\alpha$ for the case of $(-\beta, -t)$ (just
consider the pair $(x,-y)$ instead of $(x,y)$). Thus it suffices
to prove bounds of $\alpha$ for the cases when $\beta$ is
nonnegative.

\begin{lemma}\label{lemma:one}
If $\beta = 4$ and $t < 0$ then we can take $\alpha = 1$.
\end{lemma}

\begin{proof}
The hypothesis says that $x=y$ and $T(x,x) = t < 0$ so we just
take $w = x$.
\end{proof}

\begin{lemma}
\label{lemma:23rd} If $\beta = 4$ and $t > 0$, then we can take
$\alpha = 1/23$.
\end{lemma}

\begin{proof}
The hypothesis says that $x=y$ and $T(x,x) = t > 0$.  By
Lemma~\ref{lemma:ortho}, the Leech lattice contains orthogonal
minimal vectors $v_1,\dots,v_{24}$ with $v_1=x$.  Then
Lemma~\ref{lemma:orthsum}, together with $T(v_1)=t$, implies that
$T(v_i) \le -t/23$ for some $i \in \{2,\dots,24\}$.
\end{proof}

\begin{lemma}\label{lemma:25th}
If $\beta = 2$ and $t < 0$, then we can take $\alpha = 2/25$.
\end{lemma}

\begin{proof}
Let $x,y$ be minimal vectors such that $\langle x,y \rangle = 2$
and $T(x,y) = t < 0$. We know that $x-y$ is a minimal vector. Also
$$
T(x-y) = T(x,x)+T(y,y)-2T(x,y) = T(x,x)+T(y,y)-2t.
$$
Thus either $T(x,x) \leq 2t/25 = -2|t|/25$ or $T(y,y) \leq 2t/25$
or $T(x-y,x-y) \geq 4t/25-2t = -46t/25 = 46|t|/25$.  However, in
the last case we see by Lemma~\ref{lemma:23rd} that $T(v,v) \leq
-(1/23)(46|t|/25) = -2|t|/25$ for some minimal vector $v$.
\end{proof}

\begin{lemma}\label{lemma:47th}
If $\beta = 2$ and $t > 0$, then we can take $\alpha = 2/47$.
\end{lemma}

\begin{proof}
Let $x,y$ be minimal vectors such that $\langle x,y \rangle = 2$
and $T(x,y) = t > 0$. Again we have
$$
T(x-y,x-y) = T(x,x)+T(y,y)-2t.
$$
Thus either $T(x-y,x-y) \leq -2t/47$ or one of $T(x,x)$ and
$T(y,y)$ is at least $(2t-2t/47)/2 = 46t/47$.  Then by
Lemma~\ref{lemma:23rd}, we see that there exists a minimal vector
$v$ with $T(v,v) \leq (-1/23)(46t/47) = -2|t|/47$.
\end{proof}

\begin{lemma}\label{lemma:60th}
If $\beta = 0$, then we can take $\alpha =  1/60$.
\end{lemma}

\begin{proof}
By possibly exchanging $(\beta,t) = (0,t)$ with $(0,-t)$ we can
assume $t > 0$. Let $x,y$ be minimal vectors such that $\langle
x,y \rangle = 0$ and $T(x,y) = t > 0$. We have computed the
intersection numbers for the Leech lattice. The intersection
number $P_0(1/2,1/2)$ is $44$, so there exists a minimal vector
$w$ with $\langle x,w \rangle = 2$ and $\langle y,w \rangle = 2$.
We compute that
\begin{align*}
\langle x+y-w,x+y-w \rangle  &= \langle x,x\rangle +\langle
y,y\rangle+ \langle w,w \rangle +2 \langle x,y \rangle  - 2
\langle x,w \rangle  - 2 \langle y,w \rangle\\
&= 4+4+4+0-4-4\\
&= 4,
\end{align*}
so $x+y-w$ is a minimal vector. Then we compute
\begin{align*}
T(x+y-w) &= T(x,x)+T(y,y)+T(w,w)+2T(x,y) - 2T(x,w) - 2T(y,w) \\
&= T(x,x)+T(y,y)+T(w,w)+ 2t  - 2T(x,w) - 2T(y,w).
\end{align*}
If $T(x,x)$ or $T(y,y)$ or $T(w,w)$ is at most $-t/60$ we are
done. Similarly if $T(x+y-w,x+y-w) \geq 23t/60$ we are done, by
Lemma~\ref{lemma:23rd}, so we may assume none of these is true.
Then we see that $T(x,w) + T(y,w) \geq 1/2(2t-3t/60-23t/60) =
47t/60$. It follows that one of the summands is at least
$47t/120$, say $T(x,w)$ without loss of generality. But since
$\langle x,w \rangle = 2$, an application of
Lemma~\ref{lemma:47th} finishes the proof.
\end{proof}

\begin{lemma}\label{lemma:1055rd}
If $\beta = 1$ and $t > 0$, then we can take $\alpha = 4/1055$.
\end{lemma}

\begin{proof}
Let $x$ and $y$ be minimal vectors with $\langle x,y \rangle = 1$. By
Lemma~\ref{lemma:leechconfig}, there are minimal vectors $w_1$, $w_2$, and
$w_3$ satisfying $\langle x,w_i \rangle = 0$, $\langle y,w_i \rangle = 2$,
and $\langle w_i,w_j \rangle = 0$ for $i \ne j$.  It follows that $z =
2y-x-w_1-w_2-w_3$ has norm 4, so it is a minimal vector. Now,
\begin{align*}
T(z,z)
= {}& 4T(y,y)+T(x,x)+\sum_i T(w_i,w_i) - 4T(x,y) \\
&  -4\sum_i T(y,w_i) + 2\sum_i T(x,w_i) + 2\sum_{i <
j} T(w_i,w_j).
\end{align*}
Since we know $T(x,y) = t$, we get
\begin{align*} 4t = {}&
4T(y,y)+T(x,x)+\sum_i T(w_i,w_i) - T(z,z)  \\
&  - 4\sum_i T(y,w_i) + 2\sum_i T(x,w_i) + 2\sum_{i <
j} T(w_i,w_j).
\end{align*}
Next, we assume that $T(v,v) \ge -\alpha t$ for all minimal vectors $v$, and
prove that $\alpha$ cannot be less than $4/1055$.  It follows from this
assumption and Lemma~\ref{lemma:23rd} that $T(y,y) \leq 23\alpha t$.
Similarly $T(x,x)$ and $T(w_i,w_i)$ are at most $23\alpha t$, and $-T(z,z)
\leq \alpha t$ by hypothesis.  The inner products $\langle y,w_i \rangle$ are
$2$ so by Lemma~\ref{lemma:25th}, $-T(y,w_i)\le 25\alpha t/2$.  Finally,
$T(x,w_i)$ and $T(w_i,w_j)$ are at most $60 \alpha t$ by
Lemma~\ref{lemma:60th}. Therefore
$$
4t \leq 8 \cdot(23\alpha t) + \alpha t + 3\cdot4\cdot(25\alpha
t/2) + 2\cdot6\cdot(60 \alpha t) = 1055 \alpha t,
$$
and hence $\alpha \geq 4/1055$.  Thus, $T(v,v) \le -4t/1055$ for some minimal
vector $v$, since otherwise $\alpha$ could be decreased.
\end{proof}

\begin{lemma}
If $\beta = 1$ and $t < 0$, then we can take $\alpha = 4/1033$.
\end{lemma}

\begin{proof}
With the same notation as in the proof of
Lemma~\ref{lemma:1055rd} we have
\begin{align*}
T(z,z) = {}& 4T(y,y)+T(x,x)+\sum_i T(w_i,w_i) - 4T(x,y) \\
& -4\sum_i T(y,w_i) + 2\sum_i T(x,w_i) + 2\sum_{i < j}
T(w_i,w_j).
\end{align*}
Since we know $T(x,y) = t < 0$, we get as before
\begin{align*} -4t = {}&
{-4}T(y,y)-T(x,x)-\sum_i T(w_i,w_i) + T(z,z)  \\
& + 4\sum_i T(y,w_i) - 2\sum_i T(x,w_i) - 2\sum_{i <
j} T(w_i,w_j).
\end{align*}
We now use the same strategy as in Lemma~\ref{lemma:1055rd}. Namely, we
assume that for all minimal vectors $v$, $T(v,v) \ge -\alpha|t| = \alpha t$.
It follows that $-T(y,y) \leq \alpha|t| = -\alpha t$, and the same holds for
$T(x,x)$ and $T(w_i,w_i)$.  By Lemma~\ref{lemma:23rd}, $T(z,z) \leq -23\alpha
t$. This time by Lemma~\ref{lemma:47th}, $T(y,w_i)$ is at most $-47\alpha
t/2$, whereas by Lemma~\ref{lemma:60th}, $T(x,w_i)$ and $T(w_i,w_j)$ are at
least $60\alpha t$, so their negatives are at most $-60\alpha t$. Therefore
$$
-4t \leq -(23\alpha t) - 8\alpha t + 3\cdot4\cdot(-47\alpha t/2)
+ 2\cdot6\cdot(-60 \alpha t) = -1033 \alpha t,
$$
and hence $\alpha \geq 4/1033$ after cancelling $-t$ which is positive.
\end{proof}

We conclude that $\alpha = 4/1055$ satisfies the properties we
stated at the beginning of the section.

\section{The case of $E_8$}
\label{section:E8}

A very similar proof shows that the $E_8$ lattice is the unique
densest lattice packing in $\R^8$. Since the details of the proofs
are analogous, and the result was already known, we merely sketch
them in this section.

In $E_8$ there are $240$ minimal vectors of length $\sqrt{2}$, if
we normalize the length as usual so that the lattice is
unimodular. Let $\Lambda$ be a lattice of covolume $1$ that is at
least as dense as $E_8$. As in the case of the Leech lattice, we
find a suitable radial function $f$ with $r \le
\sqrt{2}(1+1.2\cdot 10^{-15})$, which proves that no sphere
packing in $\R^8$ can exceed the density of $E_8$ by a factor of
more than $1+10^{-14}$.

This function allows us to show that for $0 \neq |x| \leq
\sqrt{7}$, the length of $x$ is restricted to the set
$$
[\sqrt{2},\sqrt{2}(1+\varepsilon)) \cup (2(1-\mu),2(1+\mu)) \cup
(\sqrt{6}(1-\nu),\sqrt{6}(1+\nu)),
$$
where
\begin{align*}
\varepsilon &= 1.45\cdot 10^{-13}, \\
\mu &= 1.03 \cdot 10^{-6}, \textup{ and} \\
\nu &= 4.44 \cdot 10^{-6}.
\end{align*}
The details of this calculation are in the accompanying PARI file
\texttt{E8verifyf.txt}. All the remaining calculations from this
point on for the $E_8$ case are verified in the Maple file
\texttt{E8rest.txt}.

Define a nearly minimal vector to be a vector with length in
$[\sqrt{2},\sqrt{2}(1+\varepsilon)]$.  We form a spherical code
by rescaling the nearly minimal vectors to lie on the unit sphere
$S^7$.

We use the polynomial
$$
f_{\varepsilon}(x) = (x+1)\left(x+\frac{1}{2}\right)^2
x^2\left(x-\left(1-\frac{1}{2(1+\varepsilon)^2}\right)\right)
$$
to show as before that there are at most $240$ nearly minimal
vectors. Similarly, the analogue of Lemma~\ref{atleast196559}
goes through with a different function defined on $\R^8$, and
proves that there are exactly $240$ nearly minimal vectors.

\subsection{Spherical code}

The analogue of Section~\ref{section:sphericalcode} is that all the
inner products between the normalized nearly minimal vectors $u/|u|$
must be either $\pm 1$ or at most $6\cdot 10^{-5}$ from some element of
the set $\{-1/2,0,1/2\}$. Then further analysis gives us the following
better bounds.

For $\langle u, v \rangle \approx 1$ we have
$$
\frac{2-(1+\varepsilon)^2}{2(1+\varepsilon)^2} \le \left\langle
\frac{u}{|u|}, \frac{v}{|v|}\right\rangle \le (1+\varepsilon)^2-\frac{1}{2},
$$
whereas for  $\langle u, v \rangle \approx 0$ we have
$$
\frac{1-(1+\mu)^2}{(1+\varepsilon)^2} \le \left\langle
\frac{u}{|u|}, \frac{v}{|v|}\right\rangle \le (1+\varepsilon)^2-(1-\mu)^2.
$$
We conclude that $\sigma \le 8.89 \cdot 10^{-6}$, where $\sigma$
is the maximal error in the inner products from the spherical
code.

\subsection{Intersection numbers}

The analogue of Lemma~\ref{nearlydesign} is easily shown using
the polynomial $f_{\varepsilon}$.

\begin{lemma}\label{nearlydesignE8}
If $g \co S^{7} \to \R$ is a polynomial of total degree at most $6$,
then
$$
\left| \sum_{z \in \mathcal{C}_\Lambda} g(z) -
\frac{240}{\vol(S^{7})} \int_{S^{7}} g(z) \, dz \right| \le 3.48
\cdot 10^{-4} |g|_2,
$$
where $|g|_2$ denotes the norm on $L^2(S^{7})$.
\end{lemma}

Now, since the kissing configuration of $E_8$ is a $7$-spherical
design, we find that the intersection numbers of $E_8$ can be obtained
by solving the linear system of equations
$$
\sum_{\alpha,\beta \in S} \alpha^i \beta^j P_\gamma(\alpha,\beta) =
\frac{240}{\vol(S^{7})} \int_{S^{7}} g_{i,j}(z) \, dz.
$$
for $i,j \in \{0,1,2\}$, where $S = \{0,\pm 1/2, \pm 1\}$ is the
set of possible inner products. As before, we know the values of
$P_\gamma(\pm 1, \alpha)$ and $P_\gamma(\alpha, \pm 1)$. When we
perform the same calculation for $\Lambda$, the value of
$$
\frac{240}{\vol(S^{7})} \int_{S^{7}} g_{i,j}(z) \, dz
$$
differs by at most $30 \sigma$ from the corresponding value for $E_8$.
If we apply Lemma~\ref{lemma:approx} and compute the error introduced
into the system of equations by going from $E_8$ to $\Lambda$, we get a
bound of
$$
\left(3.48\cdot 10^{-4}\right) \frac{\pi^2}{\sqrt{3}} + 30\sigma +
240(1+2\sigma)\sigma < 4.4 \cdot 10^{-3}.
$$
The $\infty$-norm of the inverse matrix is $100$, from which it
follows that the intersection numbers in $\Lambda$ differ from
those in $E_8$ by at most $0.44$.  Because that number is less
than $1$, they must be the same.

\subsection{Association scheme}

The proof of uniqueness for the Leech lattice association scheme depended
only on the eutaxy of the Leech lattice and the uniqueness of its kissing
arrangement. The same holds for the $E_8$ lattice. With $N = 240$ and $C =
4/N = 1/60$ we have the following lemma.

\begin{lemma}
\label{lemma:eutacticE8} For every $x \in \R^{8}$,
$$
\langle x,x\rangle = C\sum_{i=1}^N \langle x,u_i\rangle^2.
$$
\end{lemma}

The proof of the lemma is essentially the same. We then use the
lemma to prove the uniqueness of the association scheme.

\begin{theorem}\label{assoc_uniquenessE8}
There is only one $4$-class association scheme with the same size, valencies
and intersection numbers as the association scheme of minimal vectors of the
$E_8$ lattice.
\end{theorem}

The proof of the theorem involves, as before, the operator
$$
P = 2C \left(A_1 + \frac{1}{2} A_{1/2} -\frac{1}{2}A_{-1/2} -
A_{-1}\right),
$$
which turns out to be a projection to an $8$-dimensional space.
Again, the proof is essentially the same as for the Leech lattice.

\subsection{Inner product bounds}

The inner product bounds use only the intersection numbers
(Lemma~\ref{lemma:quarter}, for which we needed the isomorphism of
association schemes, deals with a case that does not occur in
$E_8$). We get almost the same bounds as in
Section~\ref{section:innerbounds}. The proofs have to be slightly
modified due to the fact that the minimal norm of the $E_8$
lattice is $2$ instead of $4$ for the Leech lattice. We get the
following result:

\begin{lemma}\label{lemma:innerboundsE8}
Let $u,v$ be nearly minimal vectors.
\begin{enumerate}
\item For $\langle u, v\rangle \approx 2$ we have
$2 \leq \langle u, v\rangle  \leq 2+(5/2)\varepsilon$.
\item For $\langle u, v\rangle \approx 1$ we have
$1-(5/2)\varepsilon \leq \langle u, v\rangle \leq
1+(9/2)\varepsilon$.
\item For $\langle u, v\rangle \approx 0$ we have
$-7\varepsilon \leq \langle u, v\rangle \leq 7\varepsilon$.
\end{enumerate}
\end{lemma}

\subsection{A basis of nearly minimal vectors}

The proof that there is a basis of nearly minimal vectors is
completely analogous. Consider the basis of minimal vectors
$$
\left( \begin{matrix} 1 & 1 & 0 & 0 & 0 & 0 & 0 & 0 \\
  1 &-1 &0 &0 &0 &0 &0 &0  \\
  0 &1 &-1 &0 &0 &0 &0 &0 \\
  0 &0 &1 &-1 &0 &0 &0 &0 \\
  0 &0 &0 &1 &-1 &0 &0 &0 \\
  0 &0 &0 &0 &1 &-1 &0 &0 \\
  0 &0 &0 &0 &0 &1 &-1 &0 \\
  1/2 &1/2 &1/2 &1/2 &1/2 &1/2 &1/2 &1/2 \\
  \end{matrix} \right)
$$
of $E_8$, where the row vectors are the basis elements. Then the
corresponding vectors in $\Lambda$ (via an isomorphism of
association schemes) can be shown to be a basis for $\Lambda$.

\subsection{Computing $\alpha$}

One first proves an analogue of Lemma~\ref{lemma:ortho} by
considering the vectors $w_i = (0,\ldots,1,1,\ldots,0)$ and $v_i
= (0,\ldots,1,-1,\ldots,0)$ where the $i$ and $i+1$ coordinates
are nonzero, for $i \in \{1,3,5,7\}$.

The analogue of Lemma~\ref{lemma:orthsum} is immediate:

\begin{lemma}
\label{lemma:E8orthsum} Let $v_1,\dots,v_{8}$ be orthogonal minimal
vectors in the $E_8$ lattice.  If $\sum_{i,j} \tilde{s}_{i,j} t_{i,j} =
0$, then
$$
\sum_{i=1}^{8} T(v_i) = 0.
$$
\end{lemma}

We then proceed to use these lemmas as before to prove bounds on
$\alpha$.

\begin{lemma}\label{lemma:alphaE8}
With the same notation for $\beta$ and $t$ as before, we have the
following bounds:
\begin{enumerate}
\item If $\beta = 2$ and $t < 0$ then we can take $\alpha = 1$.
\item If $\beta = 2$ and $t > 0$ then we can take $\alpha = 1/7$.
\item If $\beta = 1$ and $t < 0$ then we can take $\alpha = 2/9$.
\item If $\beta = 1$ and $t > 0$ then we can take $\alpha = 2/15$.
\item If $\beta = 0$ then we can take $\alpha = 1/20$.
\end{enumerate}
\end{lemma}

In fact, the greatest possible value of $\alpha$ is $1/7$, and
that can be rigorously proved by computer calculations.  (The
linear programs are small enough that one can solve them using
exact rational arithmetic.  The Maple file \texttt{E8seventh.txt}
contains the calculations.) However, the weaker bound of $1/20$
will suffice.

\subsection{Local optimality of $E_8$}

The proofs of perfection and eutaxy of $E_8$ closely parallel
those for the Leech lattice. The Gram matrix for the basis of
$E_8$ that we chose above is small enough that we can compute all
its minors quickly. We find that for $0 < \rho < 10^{-3}$,
\begin{align*}
D_\rho \geq {}& 1 - \rho^2(7936 + 21162 \sqrt{3} \rho + 84256
\rho^2 + 47300 \sqrt{5} \rho^3 \\
& + 74088 \rho^4 + 10290 \sqrt{7} \rho^5 + 4096 \rho^6) \\
\geq {}& 1 - 7973 \rho^2.
\end{align*}
We conclude from the above calculations that for every
perturbation by $\rho t_{i,j}$ where $\max\{|t_{i,j}|\} \leq 1$ and
$\sum_{i,j} \tilde{s}_{i,j}t_{i,j} = 0$, we have some $k$ such that
$$
\frac{D_\rho^{-1/8} Q_\rho(u_k)}{D^{-1/8} \cdot 2} \le
(1-\rho/40) (1-7973 \rho^2)^{-1/8}.
$$
This bound is strictly less than $1$ when $0 < \rho \le 2.5
\cdot 10^{-5}$.

As before, we have one final issue to deal with: our perturbation
may not satisfy $\sum_{i,j} \tilde{s}_{i,j}t_{i,j} = 0$. We
normalize as before by setting
$$
t_{i,j} = \frac{s_{i,j} + \eta_{i,j}}{1+\Delta/8} - s_{i,j}
$$
(the notation is as in the Leech lattice case), and need to check that
$|t_{i,j}| \leq 2.5 \cdot 10^{-5}$. From
Lemma~\ref{lemma:innerboundsE8}, we have $|\eta_{i,j}| \le 7
\varepsilon$. Now the sum of the absolute values of the entries of
$\tilde{S}$ is $620$, so $|\Delta| \le 4340\varepsilon$. The maximum
value of $|s_{i,j}|$ is $2$, so putting everything together we have
$$
|t_{i,j}| = \left|\frac{\eta_{i,j} +
s_{i,j}\Delta/8}{1+\Delta/8}\right| \leq \frac{7\varepsilon +
2\cdot 4340\varepsilon/8}{1-4340\varepsilon/8} < 1.6 \cdot
10^{-10}.
$$
This proves that $\Lambda$ must be the same as $E_8$ since it
lies within the range of local optimality of $E_8$.

This completes our new proof of the optimality of $E_8$:

\begin{theorem}[Blichfeldt, Vet\v{c}inkin]\label{theorem:main2}
The $E_8$ root lattice is the unique densest lattice in $\R^{8}$,
up to scaling and isometries of $\R^{8}$.
\end{theorem}

\section*{Acknowledgements}

We thank Noam Elkies, L\'aszl\'o Lov\'asz, and Stephen D.\ Miller for
helpful discussions, Richard Pollack, Fabrice Rouillier, and
Marie-Fran\c coise Roy for advice on computer algebra, John Dunagan for
advice on linear programming software, William Stein for allowing us to
use several of his computers, Dimitar Jetchev for translating our
computer code into Magma, and Richard Borcherds, Noam Elkies, Simon
Litsyn, and Eric Rains for comments on the manuscript.

\appendix
\section{Computer calculations}
\label{appendix:computer}

The computer files for checking our calculations are available
from the \texttt{arXiv.org} e-print archive.  This paper is
available as \texttt{math.MG/0403263}.  To access the auxiliary
files, download the source files for the paper.  That will
produce not only the \LaTeX{} files for the paper but also the
computer algebra code.

By far the most extensive use of computer calculations in this paper occurs
in Subsection~\ref{howlarge}.  The calculations are carried out in the file
\texttt{verifyf.txt}, which consists of PARI code.  PARI is a free computer
algebra system designed for rapid number-theoretic calculations.  See
\url{http://pari.math.u-bordeaux.fr} for more information on PARI or to
download a copy.

Our PARI files all contain comments that should help make them understandable
to those unfamiliar with PARI.

First, we will explain how one proves the properties of the
function $f$ used in Subsection~\ref{howlarge};  then we will
explain how it was constructed.  Finally, we briefly discuss the
verification of the other calculations in this paper.

There are two auxiliary files for dealing with $f$:
\texttt{fcoeffs.txt} contains coefficients $c_0,\dots,c_{803}$ which we
will use to construct $f$, and \texttt{roots.txt} contains values
$r_0,\dots,r_{200}$ that are nearly roots (meaning the polynomial is
very near $0$ at those locations, although in fact no real roots are
nearby, except near $r_0$).

The file \texttt{verifyf.txt} carries out the following
verifications. Let
\begin{equation}\label{eq:f0}
f_0(z) = \sum_{i=0}^{803} c_i i! L_i^{11}(z).
\end{equation}
Define $f \co \R^{24} \to \R$ by $f(x) =
f_0(2\pi|x|^2)e^{-\pi|x|^2}/10^{3000}$.  (The denominator of
$10^{3000}$ makes the coefficients $c_i$ integers, which is
convenient.)  The value of $r$ used in Subsection~\ref{howlarge}
satisfies $2\pi r^2 = r_0$.

We first need to check that $f(x) \le 0$ for $|x| \ge r$, which is equivalent
to $f_0(z) \le 0$ for $z \ge r_0$.  In principle one could check this
straightforwardly using Sturm's theorem, but that takes a tremendous amount
of time for such a huge polynomial. Instead, we will use Descartes' rule of
signs in a somewhat complicated way that may not appear \emph{a priori}
superior, but works overwhelmingly better in practice: the number of roots of
a polynomial $p$ in the interval $(a,b)$ is at most the number of sign
changes in the coefficients of
$$
p\left(\frac{a+bz}{1+z}\right)(1+z)^{\deg(p)},
$$
and is congruent to it modulo $2$. For $(a,\infty)$ one can simply
use $p(a+z)$.  This result is sometimes known as Jacobi's rule of
signs (see Corollary 10.1.13 in \cite[p.~320]{RS}).

We check using Jacobi's rule of signs that $f_0$ has no roots in
$(\lceil r_{200} \rceil, \infty)$ or $(\lceil r_i \rceil, \lfloor
r_{i+1} \rfloor)$ for $1 \le i \le 199$ (and also check that it
does not vanish at the endpoints), and that it has exactly one
root in $(0,\lfloor r_1 \rfloor)$. Then we must check that it
does not vanish on $[\lfloor r_i \rfloor, \lceil r_i \rceil]$
with $1 \le i \le 200$, and that its one root is less than $r_0$
(for that we simply compute the sign of $f_0(r_0)$).

Dealing with the intervals $[\lfloor r_i \rfloor, \lceil r_i
\rceil]$ is a little more difficult.  We use Jacobi's rule of
signs to check that $f''_0$ has no roots on these intervals, so
$f'_0$ is monotonic and has at most one root in each.  We then
check that $f'_0(r_i-10^{-350})>0$ and $f'_0(r_i+10^{-350}) < 0$,
so the maximum of $f_0$ must occur within $10^{-350}$ of $r_i$.
However,
$$
|f'_0(r_i \pm 10^{-350})| \le 10^{3285},
$$
from which it follows (by the mean value theorem and the
monotonicity of $f_0'$) that
$$
\left| f_0(r_i) - \max_{x \in [\lfloor r_i \rfloor, \lceil r_i
\rceil]} f_0(x) \right| \le 10^{-350} \cdot 10^{3285} = 10^{2935}.
$$
Because $|f_0(r_i)| \le -10^{2945}$, $f_0$ has no roots in
$[\lfloor r_i \rfloor, \lceil r_i \rceil]$.

To make the proof more efficient, we arrange the calculations to
ensure that only integer arithmetic is used, so that PARI does
not spend time reducing fractions to lowest terms. (That explains
why the coefficients in \eqref{eq:f0} are multiplied by $i!$.)

Dealing with $\widehat{f}$ is similar, but of course it uses the
polynomial
$$
h_0(z) = \sum_{i=0}^{803} (-1)^i c_i i! L_i^{11}(z)
$$
instead of $f_0$, and in this case $r_0$ plays the same role as
$r_1,\dots,r_{200}$ do (rather than being treated differently, as
above).

All that remains is to check $f(0)=\widehat{f}(0)=1$, which is
true because $f_0(0) = h_0(0) = 10^{3000}$. These calculations
suffice for the proof that no sphere packing can be more than a
factor of $1+1.65\cdot10^{-30}$ times denser than
$\Lambda_{24}$.  However, Proposition~\ref{lengthbounds} requires
more detailed information on the values of $f$.  In
\texttt{verifyf.txt} we check these inequalities using
$$
\sum_{i=0}^{351} \frac{(-z)^i}{i!} \le e^{-z} \le
\sum_{i=0}^{350} \frac{(-z)^i}{i!}
$$
for $0 \le z \le 60$, as well as rational upper and lower bounds
for $\pi$, to avoid having to deal with irrational numbers.

Unfortunately the file \texttt{verifyf.txt} cannot address where $f$
comes from (it takes far longer to locate $f$ than to verify its
properties).  The construction is based on the numerical technique
described in Section~7 of \cite{CE}, which describes functions via
forced double root locations that are then repeatedly perturbed until
they reach a local optimum.  Instead of this straightforward
optimization algorithm we implemented a high-dimensional version of
Newton's method (which locates a root of the derivative of $f(2)$ as an
implicit function of the forced double roots---this is valuable because
$f(2)$ is roughly proportional to $|r-2|$). We arrived at $200$ forced
double roots, whose locations are specified in \texttt{roots.txt}.  The
next step of the method from \cite{CE} also caused trouble: solving for
the polynomial with those forced roots.  Exact rational arithmetic was
immensely time-consuming and produced huge denominators. Instead, we
carried out high-precision floating point arithmetic and rounded the
result to within $10^{-3000}$. That could ruin the sign conditions on
$f$ and $\widehat{f}$, by turning a double root into a pair of nearby
single roots, so instead of solving for actual double roots we required
that $f$ and $\widehat{f}$ should stay slightly on the correct side of
$0$.  That does not greatly change the resulting bounds, and leads to
the function $f$ used in this paper.

Lemma~\ref{atleast196559} requires the second largest amount of
computation to check, although far less than
Subsection~\ref{howlarge}. The file \texttt{verifyg.txt} and the
auxiliary file \texttt{gcoeffs.txt} deal with it in a fairly
straightforward way.

The calculations required for \eqref{eq:drho} and
\eqref{eq:upperbd} are dealt with in \texttt{verifygram.txt}.

Finally, the file \texttt{verifyrest.txt} verifies all remaining calculations
in the paper for the Leech lattice case.  Many of the calculations in this
file could be checked by hand, but that would be unpleasant and error-prone.
Instead of PARI code, this file contains Maple code.  Maple is less efficient
but more flexible, and it seemed easiest to use in these calculations.

Our computer calculations are completely rigorous, in the sense that they are
carried out using exact arithmetic, and thus avoid traps such as round-off
error.  Nevertheless, we cannot eliminate the possibility of a hardware error
or a bug in software beyond our control, such as the operating system or
computer algebra program.  Under these circumstances, even assuming our
programs are correctly written they may return incorrect results.  This
possibility appears quite unlikely, and we do not consider it a serious
worry.  However, we have addressed it by asking Dimitar Jetchev to translate
our programs into the Magma computer algebra system.  His translation is
contained in the file \texttt{magmacode.txt}.  Using it, we have
independently checked our calculations on a different type of processor, a
different computer algebra system, and a different operating system.

We have also documented our calculations for the $E_8$ proof. The
details can be found in the files \texttt{E8verifyf.txt} and
\texttt{E8rest.txt}.  We include fewer comments in these files
than in the others, because their structure is parallel to the
Leech lattice case. The file \texttt{E8seventh.txt} contains a
proof that $\alpha=1/7$ works in the $E_8$ case, although we do
not require that for the proof of Theorem~\ref{theorem:main2}.

\section{Background}
\label{appendix:background}

In this appendix we collect brief definitions and descriptions of
a few of the principal objects and techniques used in this paper.

An \textit{even unimodular lattice\/} $\Lambda \subset \R^n$ is a
lattice such that $|\Lambda|=1$, $\langle x,y \rangle \in \Z$ for all
$x,y \in \Lambda$, and $\langle x,x \rangle \in 2\Z$ for all $x \in
\Lambda$.  Such lattices exist only when $n$ is a multiple of $8$.  Up
to isometries of $\R^n$, the unique example when $n=8$ is $E_8$, and
there are $24$ examples in $\R^{24}$, among which the Leech lattice is
the unique one containing no vectors of length $\sqrt{2}$.  See
\cite[p.~48]{CS} for more information.

A \textit{spherical code\/} of minimal angle $\varphi$ in the unit
sphere $S^{n-1}$ is a collection $\mathcal{C} \subset S^{n-1}$ of
points such that $\langle x,y \rangle \le \cos \varphi$ for all
$x,y \in \mathcal{C}$ with $x \ne y$.  In other words, no two
distinct points of $\mathcal{C}$ form an angle smaller than
$\varphi$ centered at the origin.  Spherical codes are to
$S^{n-1}$ as binary error-correcting codes are to $\{0,1\}^n$, or
as sphere packings are to $\R^n$.

The most important approach to bounding the size of spherical codes is
\textit{linear programming bounds\/} (due to Delsarte \cite{D}; see
Chapter~9 of \cite{CS} for an exposition). These bounds rely on the
following property of the ultraspherical polynomials, which follows
from Lemma~\ref{lemma:nonneg}: if $\lambda = n/2-1$, then for each
finite subset $\mathcal{C} \subset S^{n-1}$,
$$
\sum_{x,y \in \mathcal{C}} C_i^\lambda(\langle x,y \rangle) \ge 0.
$$
Suppose $f_0,\dots,f_d \ge 0$, and that $f(z) = \sum_{i=0}^d f_i
C_i^\lambda(z)$ satisfies $f(z) \le 0$ for $z  \in [-1,\cos
\varphi]$. Then every spherical code $\mathcal{C}$ in $S^{n-1}$
with minimal angle $\varphi$ has size bounded by $|\mathcal{C}|
\le f(1)/f_0$ (assuming $f_0\ne0$). The proof is simple:
\begin{align*}
|\mathcal{C}|f(1) & =  \sum_{x \in \mathcal{C}} f(\langle x,x
\rangle)\\
& \ge \sum_{x,y \in \mathcal{C}} f(\langle x,y
\rangle)\\
& =  \sum_{x,y \in \mathcal{C}} f_0 + \sum_{i=1}^d f_i \sum_{x,y
\in \mathcal{C}} C_i^\lambda(\langle x,y \rangle)\\
& \ge  \sum_{x,y \in \mathcal{C}} f_0\\
& =  |\mathcal{C}|^2 f_0.
\end{align*}

The most dramatic application of the linear programming bounds for
spherical codes is the solution of the \textit{kissing problem\/} in
$\R^8$ and $\R^{24}$: how many unit balls can be placed tangent to a
central unit ball so that their interiors do not overlap?  This
condition amounts to saying that the points of tangency form a
spherical code with minimal angle $\pi/3$. In $\R^8$ the answer is
$240$, and in $\R^{24}$ the answer is $196560$. The codes are formed
from the minimal vectors in $E_8$ and the Leech lattice, respectively
(of course the radius of the sphere involved differs in the two
examples). Optimality follows from the linear programming bounds by
taking
$$
f(z) =
(z+1)\left(z+\frac{1}{2}\right)^2z^2\left(z-\frac{1}{2}\right)
$$
in $\R^8$ and
$$
f(z) =
(z+1)\left(z+\frac{1}{2}\right)^2\left(z+\frac{1}{4}\right)^2
z^2\left(z-\frac{1}{4}\right)^2 \left(z-\frac{1}{2}\right)
$$
in $\R^{24}$.  This most remarkable fact was discovered
independently by Levenshtein \cite{L} and by Odlyzko and Sloane
\cite{OS}.

A \textit{spherical $t$-design\/} in $S^{n-1}$ is a non-empty finite
subset $\mathcal{D}$ of $S^{n-1}$ such that for every polynomial $f \co
\R^n \to \R$ of total degree at most $t$,
$$
\frac{1}{|\mathcal{D}|}\sum_{x \in \mathcal{D}} f(x) =
\frac{1}{\vol(S^{n-1})} \int_{S^{n-1}} f(x) \, dx.
$$
In other words, the average of $f$ over $\mathcal{D}$ equals its
average over the entire sphere.  The minimal vectors of the Leech
lattice form a spherical $11$-design, and those of $E_8$ form a
spherical $7$-design.

A \textit{$k$-class association scheme\/} is a set $\mathcal{A}$
together with a partition
$$
\mathcal{A}^2 = \mathcal{A}_0 \cup \dots \cup \mathcal{A}_k
$$
such that $\mathcal{A}_0 = \{(x,x) : x \in \mathcal{A}\}$, $(x,y) \in
\mathcal{A}_i$ if and only if $(y,x) \in \mathcal{A}_i$, and the following
property holds. Fix $\ell$, $m$, and $n$ in $\Z \cap [0,k]$; then for all
$x,y \in \mathcal{A}$ with $(x,y) \in \mathcal{A}_\ell$, there are the same
number $P_\ell(m,n)$ of $z \in \mathcal{A}$ such that $(x,z) \in
\mathcal{A}_m$ and $(y,z) \in \mathcal{A}_n$. (That is, $P_\ell(m,n)$ depends
only on $\ell$, $m$, and $n$, and not on $x$ and $y$.)  These numbers are
called the \textit{intersection numbers\/} of the association scheme.  When
$\ell=0$ and $m=n$ they are also called \textit{valencies\/}.  In the body of
the paper we modify this notation slightly for the association scheme
${\mathcal C}_\Lambda$ (in a trivial way): we label the classes with the
corresponding inner products in $\mathcal{C}_{24}$, and use these labels in
the notation for the intersection numbers.

\end{document}